\documentclass[]{article}
\catcode`\@=11
\catcode`\@=12
\usepackage{amsmath}
\usepackage{amssymb}
\usepackage{amsthm}

\usepackage{amsmath,amsfonts,amssymb,graphicx, mathrsfs}

\usepackage{color,graphicx}

\newtheorem{theorem}{Theorem}[section]

\newtheorem{condition}[theorem]{Condition}

\newtheorem{lemma}[theorem]{Lemma}

\newtheorem{proposition}[theorem]{Proposition}
\newtheorem{remark}[theorem]{Remark}

\begin{document}

\title{Weak Solutions of $3D$ Compressible Navier-Stokes Equations in Critical Case}
\author{P.I. PLotnikov\\}
%{\small Lavrentyev Institute of Hydrodynamics  \\
 %             \email: piplotnikov@mail.ru}}
\date{}

\maketitle

\tableofcontents

\begin{abstract}
New estimates  of the potentials  of solutions  to the compressible
Navier-Stokes equations are derived. The result obtained are applied
to  boundary value problems for the compressible Navier-Stokes
equations with the critical adiabatic exponents. The cancelation of
concentrations of the kinetic energy density is proved.
\end{abstract}

\section{Introduction}\label{section1}
\renewcommand{\theequation}{1.\arabic{equation}}
\setcounter{equation}{0} Suppose a viscous compressible fluid
occupies
  a bounded domain $\Omega\subset \mathbb{R}^3$,
The  state of the fluid is characterized  by the macroscopic
quantities: the {density}
 $\varrho(x,t)$ and the {velocity}
 $\mathbf{u}(x,t)$.
  The problem is to  find
$\mathbf u(x,t)$ and $\varrho(x,t)$ satisfying the following
equations and boundary conditions in the cylinder $Q_T=\Omega\times
(0,T)$.

 \begin{subequations}
  \label{stokes1}
\begin{gather}
\label{stokes2} \partial_t (\varrho \mathbf{u})+\text{~div~}(\varrho
\mathbf u \otimes\mathbf{u})+\nabla
p(\varrho)=\text{~ div~}\mathbb S(\mathbf u)+\varrho\mathbf f\quad\text{in~~}Q_T,\\
 \label{stokes3}\partial_t\varrho+\text{~ div~}(\varrho \mathbf{u})=0,\quad \varrho\geq
 0\text{~~in~~}Q_T,\\\label{stokes4}
\mathbf u=0\quad\text{on~} \partial\Omega\times (0,T),\\
\label{stokes5} \mathbf u(x,0)=\mathbf u_0(x),\quad
\varrho(x,0)=\varrho_0(x)> 0\quad\text{in~} \Omega.
\end{gather}
Here, the vector field $\mathbf f\in L^\infty (Q_T)$  denotes the
density of external mass forces,
 the {viscous stress
tensor}
 $\mathbb S(\mathbf u)$ has the form
\begin{equation}
\label{stokes4a}
 \mathbb S(\mathbf{u})=\nu_1 \big(\nabla \mathbf{u}+\nabla\mathbf{u}^\top\big)+ \nu_2 \text{div~}
 \mathbf u\,\mathbb{I},
\end{equation}
\end{subequations}
where the viscosity coefficients satisfy the inequalities $\nu_1>0$
and $\nu_1+\nu_2\geq 0$.
 The first nonlocal results
concerning the mathematical theory of compressible Navier-Stokes
equations are due to P.-L. Lions. In monograph \cite{PLL} he
established the existence of a  renormalized solution to
nonstationary boundary value problem for the Navier-Stokes equations
with the pressure function $p\sim\varrho^\gamma$ for all
$\gamma>5/3$.
 More recently, Feireisl, Novotn\'{y},  and Petzeltov\'{a}, see
  \cite{FNP},  proved the existence result for all
 $\gamma>3/2$,
 see also monographs \cite{FEIRBOOK}, \cite{NovotnyStra}
 for references and details. For   $\gamma\leq 3/2$,
The question on solvability of problem \eqref{stokes1} with critical
and subcritical $\gamma$ is still open. In this case the problem of
solution existence was posed in \cite{AAL}. For two-dimensional
flows, this problem was studied in \cite{padula1} and \cite{PLWe}.
Accordingly, in what follows, we consider three-dimensional flows.

The main difficulty is the so called concentration problem, see
\cite{PLL} ch.6.6 and \cite{AAL}. It is easy to see that in the
three-dimensional case the
 energy estimates and embedding theorems guarantee the inclusion
$\varrho|\mathbf u|^2\in L^s(\Omega)$ with $s>1$ if and  only if
$\gamma>3/2$. Hence, for $\gamma\leq 3/2$ we have only an $L^1$
estimate for the density of the kinetic energy.

The problem can be formulated as follows.  Choose the artificial
pressure function in the form
\begin{equation}\label{pressure1}
p(\varrho)=\varrho^{\gamma} P_\varepsilon (\varrho)\equiv
\varrho^{\gamma} (1+\varepsilon\varrho^m), \quad m>1,
\end{equation}
where $\varepsilon>0$ is a small parameter. The question is: Under
what conditions is a weak limit of solutions $\varrho_\varepsilon,
\mathbf u_\varepsilon$ to the regularized equations with
$p=\rho^\gamma P_\varepsilon$ a solution to problem \eqref{stokes1}
with $p=\rho^\gamma$. It is known, \cite{FNP}, that the sequence of
solutions $\varrho_\varepsilon,\mathbf u_\varepsilon$ to the
regularized problem satisfies the energy inequality
\begin{equation*}
    \|\varrho_\varepsilon^{\gamma}\|_{L^\infty(0,T; L^1(\Omega))}+
 \|\varrho_\varepsilon |\mathbf u_\varepsilon|^2\|_{L^\infty(0,T; L^1(\Omega))}
 + \|\mathbf u_\varepsilon\|_{L^2(0,T; W^{1,2}(\Omega))}\leq c,
\end{equation*}
where $c$ is independent of $\varepsilon$. We may assume, after
passing to a subsequence, that
\begin{gather*}
 \mathbf u_\varepsilon\to\mathbf u\text{~~weakly in~~}L^2(0,T;
W^{1,2}_0(\Omega)), \\\varrho_\varepsilon\to\varrho\text{~~star
weakly in~~}L^\infty(0,T; L^{\gamma}(\Omega)),
\\
\varrho_\varepsilon \mathbf u_\varepsilon\otimes\mathbf
u_\varepsilon\to \mathcal{M} \text{~~star weakly in the space of
Radon measures as~~} \varepsilon \to 0,
\end{gather*}
where $\mathcal{M}=(\mathcal{M}_{ij})_{3\times 3}$ denotes a
matrix-valued Radon measure in $\Omega$. In the general case the
weak star defect measure $\mathcal{M}- \varrho\mathbf u\otimes
\mathbf u\neq 0$. This  leads to the so-called concentration
problem, which was widely discussed in the mathematical literature
in connection with vortex sheets dynamics.
 Hence the question is to describe
the structure of the defect measure and to find conditions under
which it  equals  $0$.

Recently, the Hausdorff and parabolic dimensions of the support of
the defect measure were estimated in \cite{Xi} and \cite{Xi1}.
Specifically, it was shown that, in the critical case, the kinetic
energy tensor is concentrated on a set of Hausdorff dimension of at
most 3/2. We investigate this issue for a sequence of kinetic energy
tensors in the case of the critical adiabatic exponent.
 Our goal is to investigate the concentration
problem in the critical case $\gamma=3/2$ and to prove the
cancelation of concentrations. It suffices to show that there is a
continuous function $\Psi:\mathbb R^+\to \mathbb R^+$ such that
$$
\Phi(s)/s\to \infty\text{~~as~~} s\to \infty
$$
and
$$
\int_{K}\Psi(\varrho_\varepsilon|\mathbf u_\varepsilon|^2)\,
dxdt\leq c(K),
$$
for every compact set $K\subset Q_T$. Here $c(K)$ is independent of
$\varepsilon$. We prove that the kinetic energy densities of weak
solutions to regularized problem \eqref{stokes1} with $\gamma=3/2$
and $p=\rho^{3/2} P_\varepsilon$ are uniformly bounded in some
logarithmic Orlitz space. Further, we will assume that the flow
domain and the given data satisfy the following condition.
\begin{condition}
\label{stokes6} The flow domain  $\Omega\subset \mathbb{R}^3$ is a
bounded domain with $C^2$ boundary.  The  data $ \varrho_0,
\mathbf{u}_0\in L^\infty(\Omega)$ and $\mathbf f\in L^\infty(Q_T)$
admit the estimate
\begin{equation}\label{stokes8}
\|\mathbf{u}_0\|_{W^{1,2}_0(\Omega)}+
\|\varrho_0\|_{L^\infty(\Omega)} +\|\mathbf f\|_{L^\infty(Q_T)}\leq
c_e,\quad \varrho_0>c>0,
\end{equation}
where $c_e$ and $ c$ are  positive constants.
\end{condition}
Hereinafter we will denote  by $c$    generic constants depending
only on $ \Omega, T, \|\varrho_0\|_{L^\infty(\Omega)}$,   $\|\mathbf
u_0\|_{L^2(\Omega)}, $  $ \|\mathbf f\|_{L^\infty(Q_T)}$, and
$\nu_i$.

Recall that  functions $\varrho$, $\mathbf u$ satisfy equations
\eqref{stokes1} in the weak sense if the integral identities
\begin{multline*}
\int_{Q_T}\big(\varrho \mathbf{u}
\cdot\partial_t\boldsymbol{\xi}\big)\, dxdt + \int_{Q_T}\big(\varrho
\mathbf u\otimes\mathbf
u+\varrho^{3/2}P_\varepsilon(\varrho)\,\mathbb I- \mathbb
S(\mathbf{u})
 \big):\nabla
\boldsymbol{\xi}\, dxdt+\\\int_{Q_T}\varrho\mathbf{f}\cdot
\boldsymbol{\xi}\, dxdt  +\int_\Omega\varrho_0(x)
\mathbf{u}_0(x)\cdot\boldsymbol{\xi}(x,0)\, dx=0
\end{multline*}
\begin{equation*}
\int_{Q_T}\Big(\varrho\partial_t\psi+ \big(\varrho \mathbf u\big)
\cdot\nabla\psi\Big)\,dxdt +\int_{\Omega }( \psi\varrho_0)(x,0)\,
dx= 0,
\end{equation*}
hold true for all vector fields
 $\boldsymbol{\xi}\in C^\infty(Q_T)$ and functions  $\psi\in C^\infty(Q_T)$,
 satisfying conditions
$ \boldsymbol{\xi}(x,T)=0$, $\psi(X,T)=0$ in $\Omega$  and
 $\boldsymbol{\xi}=0$  on  $ \partial\Omega\times (0,T)$.

The main goal of this paper is the proof of the following
\begin{theorem}\label{stokesthm}
Let $0<\sigma<8/21$ and a compact set $K\subset Q_T$. Furthermore
assume  that a weak solution $\varrho_\varepsilon$, $\mathbf
u_\varepsilon$ to problem \eqref{stokes1} with $p(\rho)=\rho^{3/2}
P_\varepsilon$ admits the energy estimate
\begin{equation*}
    \|\varrho_\varepsilon^{3/2}\|_{L^\infty(0,T; L^1(\Omega))}+
 \|\varrho_\varepsilon |\mathbf u_\varepsilon|^2\|_{L^\infty(0,T; L^1(\Omega))}
 + \|\mathbf u_\varepsilon\|_{L^2(0,T; W^{1,2}(\Omega))}\leq E
\end{equation*}
Then the kinetic energy density $\varrho_\varepsilon|\mathbf
u_\varepsilon|^2$ satisfies the equality
\begin{equation}\label{stokes5abcd}
    \int_{K}\varrho_\varepsilon|\mathbf u_\varepsilon|^2
    \big(\ln (1+\varrho_\varepsilon|\mathbf u_\varepsilon|^2 )\big)^{\sigma}\leq c,
\end{equation}
where $c$ depends only on $E$, $K$, $\nu_i$,  and $Q_T$.
\end{theorem}
The rest of the paper is devoted to the proof of this theorem.
\section{Localization. Auxiliary problem}\label{pont}
\renewcommand{\theequation}{2.\arabic{equation}}
\setcounter{equation}{0}

Since estimate \eqref{stokes5abcd} is local, it is convenient to
reduce the problem to the case when a solution is defined in the
whole space $\mathbb R^4$. In our studies, the initial and boundary
conditions do not play any role.  Further,  we will consider  weak
solutions $\varrho_\varepsilon,\mathbf u_\varepsilon:Q_T\to
R^+\times R^3$ to equations \eqref{stokes2}-\eqref{stokes3} with the
pressure function
\begin{equation}\label{pressure2}
p(\varrho)=\varrho^{3/2} {P}_\varepsilon(\varrho), \quad
{P}_\varepsilon(\varrho)=1+\varepsilon \varrho^m, \quad m>1.
\end{equation}
We  assume that these solutions admit the estimate
\begin{equation}\label{pont1}
    \|\varrho_\varepsilon^{3/2}{P}_\varepsilon(\varrho_\varepsilon)\|_{L^\infty(0,T;
    L^1(\Omega))}+
 \|\varrho_\varepsilon |\mathbf u_\varepsilon|^2\|_{L^\infty(0,T; L^1(\Omega))}
 + \|\mathbf u_\varepsilon\|_{L^2(0,T; W^{1,2}(\Omega))}\leq c.
\end{equation}
It is also assumed that the integral identities
\begin{multline}
\int_{Q_T}\varrho_\varepsilon \mathbf{u_\varepsilon}
\cdot\partial_t\boldsymbol{\xi}\, dxdt +
\int_{Q_T}\big(\varrho_\varepsilon \mathbf
u_\varepsilon\otimes\mathbf
u_\varepsilon+\varrho_\varepsilon^{3/2}{P}_\varepsilon(\varrho_\varepsilon)\,\mathbb
I- \mathbb S
 \big):\nabla
\boldsymbol{\xi}\,
dxdt\\+\int_{Q_T}\varrho_\varepsilon\mathbf{f}\cdot
\boldsymbol{\xi}\, dxdt=0 \label{navier8in}
\end{multline}
\begin{equation}\label{navier10in}
\int_{Q_T}\Big(\varrho_\varepsilon\partial_t\psi+
\big(\varrho_\varepsilon \mathbf u_\varepsilon\big)
\cdot\nabla\psi\Big)\,dxdt = 0,
\end{equation}
 hold true for all vector fields
 $\boldsymbol{\xi}\in C^\infty_0(Q_T)$ and functions  $\psi\in
 C^\infty_0(Q_T)$.
\begin{remark}\label{S}
Notice that  the specific form of the viscous stress tensor $\mathbb
S$ is not essential. It is suffices to assume that it is defined in
$Q_T$ and admits the estimate
\begin{equation}\label{S2}
    \|\mathbb S\|_{L^2(Q_T)}\leq c.
\end{equation}
\end{remark}
In order to localize the problem,  fix an arbitrary compact set
$K\subset Q_T$. Next, fix an arbitrary nonnegative function
$\lambda\in C^\infty_0(\mathbb R^4)$ such that
$$
\text{supp~}\lambda\subset Q_T\text{~~and~~} \lambda(x,t)=1
\text{~~in a neighborhood of ~~} K.
$$
Now set
$$
\rho(x,t)=\lambda^4(x,t)\, \varrho_\varepsilon(x,t), \quad \mathbf
v(x,t)=\lambda(x,t)\, \mathbf u_\varepsilon(x,t).
$$
Extend  $(\rho, \mathbf v)$  and   $(\varrho_\varepsilon,\mathbf
u_\varepsilon)$ by zero to the whole space $ R^4$. The following
lemma constitutes the momentum balance equation for the functions
$\rho$ and $\mathbf v$.
\begin{lemma}\label{mambda1}
Under the above assumptions, the functions $\rho, \mathbf v$
satisfies the integral identity
\begin{multline}
\int_{R^4}\lambda\rho \mathbf{v} \cdot\partial_t\boldsymbol{\xi}\,
dxdt + \int_{R^4}\big(\rho \mathbf v\otimes\mathbf
v+\rho^{3/2}{P}_\varepsilon\,\mathbb I
 \big):\nabla
\boldsymbol{\xi}\, dxdt\\
-\int_{R^4} \mathbb T :\nabla \boldsymbol{\xi}\,
dxdt+\int_{R^4}\mathbf{F}\cdot \boldsymbol{\xi}\, dxdt=0
\label{mambda2}
\end{multline}
for all $\boldsymbol\xi\in C^\infty_0(R^4)$. The functions $\rho$,
$\mathbf v$, $\mathbb T$, $\mathbf F$, $P_\varepsilon$ are compactly
supported in $Q_T$ and admit the estimates
\begin{gather}\label{mambda3}
  \|\rho^{3/2}{P}_\varepsilon\|_{L^\infty(R;
    L^1(R^3))}\leq c, \quad
    P_\varepsilon=(1+\varepsilon\varrho^m),\\\label{mambda4}
 \|\rho |\mathbf v|^2\|_{L^\infty(R; L^1(R^3))}
 + \|\mathbf v\|_{L^2(R; W^{1,2}(R^3))}\leq c,\\\label{mambda5}
\|\mathbb T\|_{L^2(R^4)}+\|\mathbf F\|_{L^1(R^4)}\leq c
\end{gather}
where the  constant $c$ is independent of $\varepsilon$.
\end{lemma}

\begin{proof}
Replace the functions $\boldsymbol \xi$ in integral identity
\eqref{navier8in} by the functions $\lambda^6\boldsymbol\xi$.
Straightforward calculations give    integral identity
\eqref{mambda2} with
\begin{gather*}
\mathbb T=\lambda^6\mathbb S, \quad \mathbf F=\mathbf F_1+\mathbf
F_2,\quad \mathbf F_1=-6\lambda^5\mathbb S\nabla\lambda,\\
\mathbf F_2=6\lambda^5\partial_t\lambda\, \varrho_\varepsilon\mathbf
u_\varepsilon +6\lambda^5(\varrho_\varepsilon\mathbf
u_\varepsilon\otimes\mathbf
u_\varepsilon)\,\nabla\lambda+6\lambda^5\varrho_\varepsilon
P_\varepsilon \nabla\lambda+\lambda^6\varrho_\varepsilon \mathbf f.
\end{gather*}
It remains to note that energy estimate \eqref{pont1} yields
\begin{equation}\label{pont6}
\|\mathbf F_1\|_{L^2(\mathbf R^4)}+ \|\mathbf
F_2\|_{L^\infty(\mathbf R; L^1(\mathbb R^3))} \leq C,
\end{equation}
where $C$  depends only on $E$, $K$, $\nu_i$, and $Q_T$.
\end{proof}
Hence Theorem \ref{stokesthm} can  be reformulated as follows
\begin{theorem}\label{pont7}
Under the assumptions of Theorem \ref{stokesthm}, the modified
energy density $\rho|\mathbf v|^2$ admits the estimate
\begin{equation}\label{stokes5abc}
    \int_{\mathbb R^4}\rho|\mathbf v|^2 \big(\ln (1+\rho|\mathbf v|^2 )\big)^{\sigma}\leq c,
\end{equation}
where $c$ depends only on $E$, $K$, $\nu_i$,  and $Q_T$.
\end{theorem}

\section{Integral identities }\label{slava}
\renewcommand{\theequation}{3.\arabic{equation}}
\setcounter{equation}{0}

In this section we derive the special integral identities for
solutions to problem  \eqref{stokes1}. We also use these identities
in order to obtain preliminary estimates for solutions to this
problem. To this end, we introduce projections, depending on a
vector $\eta\in \mathbb R^3\setminus\{0\}$,
\begin{equation}\label{slava1}
    \Pi^\bot(\eta)x\,=\, x\,-\,\frac{\eta\cdot x}{|\eta|}\,\eta, \quad
    \Pi^\parallel(\eta)x\,=\,\frac{\eta\cdot x}{|\eta|}\,\eta.
\end{equation}
Let us consider the functions $\rho, \mathbf v$ satisfying
conditions of Lemma \ref{mambda1}. Decompose the difference $\mathbf
v(x,t)-\mathbf v(y,t) $ into two orthogonal vectors
$$
\mathbf v(x,t)-\mathbf v(y,t)=\Pi^\parallel(x-y)(\mathbf
v(x,t)-\mathbf v(y,t))+\Pi^\bot(x-y)(\mathbf v(x,t)-\mathbf v(y,t))
$$
and  introduce the "energy functions"

\begin{equation}\label{mila1}\begin{split}
\mathbf
E^\parallel(x,y,t)=\frac{1}{2}\Big(\,\Pi^\parallel(x-y)(\mathbf
v(x,t)-\mathbf v(y,t))\Big)^2+\varrho(x,t)^{1/2}P_\varepsilon(x,t), \\
\mathbf E^\bot(x,y,t)=\frac{1}{2}\Big(\,\Pi^\bot(x-y)(\mathbf
v(x,t)-\mathbf
v(y,t))\Big)^2+2\varrho(x,t)^{1/2}P_\varepsilon(x,t),\\
\mathbf E=\mathbf E^\bot+\mathbf E^\parallel=\frac{1}{2}(\mathbf
v(x,t)-\mathbf v(y,t))^2+3\varrho(x,t)^{1/2}P_\varepsilon(x,t).
\end{split}\end{equation}
Finally, introduce the functions $\Gamma, W:(0,e^{-1})\to R$ defined
by the equalities
\begin{equation}\label{slava9}
    \Gamma(r)\,=\,\frac{1}{r}\int_{R^4}\rho(x,t)\,
    \Big(\int_{B(x,r)}\rho(y,t)\,\mathbf E^\parallel(x,y,t)\,dy\Big)dxdt,
\end{equation}
\begin{equation}\label{slava10}
    W(r)\,=\,\int_0^r\Big\{\int_{R^4}\rho(x,t)\,
    \Big(\int_{B(x,s)}\rho(y,t)\,\mathbf
    E^\bot(x,y,t)\,dy\Big)dxdt\Big\}\frac{ds}{s^2}.
\end{equation}
The following proposition is the main result of this section.
\begin{proposition}\label{slava14p}
Under the assumptions of Lemma \ref{mambda1}, the functions $W$,
$\Gamma$ and
\begin{equation}\label{slava30}
    \Phi=\Gamma-W
\end{equation}
admit the estimates
\begin{equation}\label{slava14}
    W(r)\,+\,\Gamma(r)\,\leq\,C\, , \,\,r\in[0,e^{-1}]
\end{equation}
\begin{equation}\label{slava15}
    |\Phi(r)|\,\leq\,c\,\sqrt{r}\, , \,\,r\in[0,e^{-1}],
\end{equation}
where $c$ is independent of $\rho$, $\mathbf v$ and $r$.
\end{proposition}
The rest of the section is devoted to the proof of Proposition
\ref{slava14p}. The proof falls into four steps.
\paragraph{Step 1. Kernels and test functions}
Fix an arbitrary $\mu\in (0,1]$ and $r>0$.  Introduce the vector
valued kernels
\begin{equation}\label{slava50}\begin{split}
K_0(x)=\Big(\frac{1}{|x|^\mu}-\frac{1}{r^\mu}\Big)\, x\text{~~for~~}
0<|x|\leq r, \quad K_0(x)=0\text{~~for~~}|x|>r,\\
K_j(x)=\partial_j K_0(X)\text{~~for~~} 0<|x|\leq r, \quad
K_j(x)=0\text{~~for~~} |x|>r. \end{split}\end{equation} Without loss
of generality we may assume that $K_0(0)=K_j(0)=0.$ The components
$K_{ij}$ of the vector valued functions $K_j$ form symmetric
$(3\times 3)$-matrix with entries
\begin{equation}\label{slava51}
K_{ij}(x)=\frac{|x|^2\delta_{ij}-\mu x_ix_j}{|x|^{2+\mu}}-
\frac{\delta_{ij}}{r^\mu}\text{~~for~~}|x|\leq r,\,\,\, K_{ij}(x)=0
\text{~~for~~} |x|>r.
\end{equation}
Introduce the vector field
\begin{equation}\label{slava33}
    \boldsymbol\xi(x,t)=\int_{
    R^3}K_0(x-y)\rho(y,t)\, dy\equiv K_0* \rho,\quad
    \partial_j\boldsymbol\xi=K_j*\rho.
\end{equation}
Note that for $\mu<1$, these functions are bounded, i.e.,
\begin{equation}\label{slava34}
    |\boldsymbol \xi(x,t)|\leq c, \quad
    |\partial_j\boldsymbol\xi(x,t)|\leq c,
\end{equation}
where the constant $c$ does not depend on $x,t$, but depends on
$\mu$.  Since $\rho^{3/2}$ belongs to the class $L^\infty(R;
L^1(R^3))$ and is compactly supported in the cylinder $Q_T$,  the
first estimate is obvious. The second is obviously true for $x\in
R^3\setminus \Omega$. For $x\in \Omega$ and
$R=2\text{~diam~}\Omega$, we have
\begin{multline}\label{slava34a}
|\partial_j\boldsymbol\xi|\leq c\int_{B(x,R)}
\frac{\rho(y,t)dy}{|x-y|^\mu}\leq\\
\Big(\int_{B(0,R)}\frac{dz}{|z|^{3\mu}}\Big)^{1/3}
\Big(\int_{B(x,R)}\varrho^{3/2}(y,t)dy\Big)^{2/3} \leq
c(1-\mu)^{-1/3}.
\end{multline}
In order to estimate  $\partial_t\boldsymbol \xi$, we use integral
identity \eqref{navier10in} which represents the weak form of the
mass balance equation.   Choose an arbitrary function $\eta\in
C^\infty_0(R^4)$ and set
$$
\psi=(\psi_1,\psi_2, \psi_3)=\lambda^4\, K_0*\eta\in
C^\infty_0(R^4).
$$
Substituting $\psi$ into \eqref{navier10in}, recalling that $\mathbf
u$, $\varrho$ are extended by zero outside of $Q_T$ and noting that
$\rho=\lambda^4\varrho$, $\mathbf v=\lambda \mathbf u$, we arrive at
the integral identity
\begin{equation}\label{slava35}\begin{split}
\int_{R^4}\big(\,\rho \, K_0*(\lambda\eta)_t+\rho u_{\varepsilon
j}\,
K_j*(\lambda\eta)\,\big)\, dxdt+\\
\int_{R^4}\big(\,\varrho_\varepsilon\partial_t\lambda^4+\varrho_\varepsilon
\mathbf u_\varepsilon\cdot \nabla\lambda^4\,\big) K_0(\lambda\eta)\,
dxdt=0.
\end{split}\end{equation}
It follows from the skew-symmetric property of the kernel $K_0$ and
the definition of $\boldsymbol\xi$ that
\begin{equation}\label{slava36}
    \int_{R^4} \rho K_0*(\lambda\eta)_t\, dxdt= -\int_{R^4}(\lambda\eta)_t\boldsymbol\xi\,
    dxdt.
\end{equation}
Next,  the symmetric property of the kernel $K_j$ implies
\begin{equation}\label{slava37}\begin{split}
\int_{R^4}\rho u_{\varepsilon j}\, K_j*(\lambda\eta)\,
dxdt=\int_{R^4} \lambda\eta\, K_j*(\rho u_{\varepsilon j})\, dxdt=\\
\int_{R^4}\eta K_j*(\rho v_j)\,dxdt+\int_{R^4}\eta[\lambda,
K_j](\rho u_{\varepsilon j})\, dxdt,
\end{split}\end{equation}
 where the commutator  is defined by the equality
\begin{equation*}
[\lambda, K_j]\cdot= \lambda\, K_j*(\cdot)- K_j*(\lambda\cdot).
\end{equation*}
Substituting equalities \eqref{slava36}-\eqref{slava37} into
\eqref{slava35}  we arrive at the integral identity
\begin{equation}\label{slava39}\begin{split}
\int_{R^4}\big(-\lambda\boldsymbol\xi\, \partial_t\eta+K_j*(\rho v_j)\,\big)\, dxdt+\\
\int_{R^4}\big(\lambda
K_0*(\partial_t(\lambda^4)\varrho+\varrho\mathbf
u_\varepsilon\cdot\nabla\lambda^4)+[\lambda, K_j](\rho
u_{\varepsilon j})-\lambda_t\boldsymbol\xi \,\big)\eta\, dxdt=0.
\end{split}\end{equation}
Now introduce the function
\begin{equation}\label{slava40}
    G=\lambda
K_0*(\partial_t(\lambda^4)\varrho+\varrho\mathbf
u_\varepsilon\cdot\nabla\lambda^4)+[\lambda, K_j](\rho
u_{\varepsilon j})-\lambda_t\boldsymbol\xi
\end{equation}
With this notation relation \eqref{slava33} is equivalent to
equality
\begin{equation}\label{slava41}
    \partial_t(\lambda \boldsymbol\xi)= -K_j*(\rho
    v_j)-G\text{~~in~~} R^4,
\end{equation}
which is understood in the sense of distributions. The following
lemma constitutes the integrability of the right hand side of this
equality.
\begin{lemma}\label{slaval2}Under the above assumptions, we have
\begin{equation}\label{slava42}\begin{split}
\|\varrho_\varepsilon\mathbf u_\varepsilon\|_{L^4(R;
L^1(R^3))}+\|\varrho_\varepsilon \mathbf
u_\varepsilon\|_{L^\infty(R^1; L^{6/5}(R^3))}\leq c, \\
 \|K_j*(\rho v_j)\|_{L^\infty(Q_T)}+\|G\|_{L^\infty(R^4)}\leq c.
\end{split}\end{equation}
\end{lemma}
\begin{proof}
Recall that $\varrho_\varepsilon$, $\mathbf u_\varepsilon$ and
$\rho$, $\mathbf v$ are supported in $Q_T$. Set
$$
\alpha=\frac{1}{2}, \quad \beta=\frac{3}{4}, \quad
\gamma=\frac{1}{4}, \quad r=4, \quad \tau=2, \quad s=4.
$$
The equality $\varrho_\varepsilon|\mathbf
u_\varepsilon|=\varrho_\varepsilon^\beta|\mathbf
u_\varepsilon|^\alpha(\varrho|\mathbf u_\varepsilon|^2)^\gamma$
along with the H\"{o}lder inequality implies
\begin{equation*}\begin{split}
\int_{R^3}\rho_\varepsilon(t)|\mathbf u_\varepsilon(t)|\, dx \leq
\Big(\int_{R^3}\varrho_\varepsilon^{\beta
\tau}dx\Big)^{1/\tau}\Big(\int_{R^3}|\mathbf u_\varepsilon|^{\alpha
r}dx\Big)^{1/r}\Big(\int_{R^3}(\varrho_\varepsilon|\mathbf
u_\varepsilon|^2)^{\gamma s }dx\Big)^{1/s}\\=
\Big(\int_{R^3}\varrho_\varepsilon(t)^{3/2}dx\Big)^{1/4}\Big(\int_{R^3}|\mathbf
u_\varepsilon(t)|^{2}dx\Big)^{1/4}\Big(\int_{R^3}\varrho_\varepsilon(t)|\mathbf
u_\varepsilon(t)|^2)dx\Big)^{1/4}\\\leq c\Big(\int_{R^3}|\mathbf
u_\varepsilon(t)|^{2}dx\Big)^{1/4}=\|\mathbf
u_\varepsilon\|_{L^2(\mathbb R^3)}^{1/2}.
\end{split}\end{equation*}
By  virtue of the energy estimate, we have $\|\mathbf
u_\varepsilon\|_{L^2(R^4)}\leq c$ which yields the desired estimate
for $\|\varrho_\varepsilon\mathbf u_\varepsilon\|_{L^4(R;
L^1(R^3))}$. Next we have for every $t\in R$,
\begin{equation*}\begin{split}
\int_{R^3}|\varrho_\varepsilon(t)\mathbf
u_\varepsilon(t)|^{6/5}dx\leq
\Big(\int_{R^3}\varrho_\varepsilon(t)^{3/2}\,
dx\Big)^{2/5}\Big(\int_{R^3}\varrho_\varepsilon|\mathbf
u_\varepsilon(t)|^2\, dx\Big)^{3/5}\leq c.
\end{split}\end{equation*}
In order to estimate the convolution of  $K_j$ and $\rho\mathbf v$,
note that
$$
|K_j*(\rho\mathbf v_j)(x,t)|\leq
c\int_{R^3}\frac{1}{|x-y|}\rho(y,t)|\mathbf v(y,t)|dy.
$$
Since $\rho(\cdot,t)\mathbf v(\cdot,t)$ is compactly supported in
$\Omega$, we have
$$
\|K_j*(\rho v_j)(\cdot,t)\|_{L^6(\Omega)}\leq c\|K_j*(\rho\mathbf
v_j)(\cdot,t)\|_{W^{2,6/5}(\Omega)} \leq c\|\rho(\cdot,t)\mathbf
v(\cdot,t)\|_{L^{6/5}(\Omega)}\leq c,
$$
which yields the desired estimates for the convolution. It remains
to note that the kernel
$$
(\lambda(x,t)-\lambda(y,t))K_j(x-y)
$$
of the commutator $[\lambda, K_j]$ is bounded by a constant
depending on $\lambda$ only. The desired estimate for $G$ follows
from this, the boundedness of the kernel $K_0$, the vector field
$\boldsymbol \xi$, and the integrability of functions
$\varrho_\varepsilon$, $\varrho_\varepsilon\mathbf u_\varepsilon$.
\end{proof}

It follows from  estimates \eqref{slava34}, representation
\eqref{slava41} for $\partial_t\boldsymbol\xi$, and Lemma
\ref{slaval2} that all integrals in integral identity
\eqref{mambda2} are well defined  for the test vector field
$\boldsymbol \xi$ given by equality \eqref{slava33}. Substituting
$\boldsymbol\xi$ and representation \eqref{slava41} into
\eqref{mambda2}  we arrive at the equality
\begin{equation}\label{slava42a}
    \mathbf I_1+\mathbf I_2+\Phi=0,
\end{equation}
where
\begin{equation}\label{slava43}\begin{split}
    \mathbf I_1=\int_{R^4}\big(\rho(\mathbf v\otimes\mathbf
    v):\nabla\boldsymbol\xi-\varrho\mathbf v \cdot K_j*(\rho
    v_j)\big)\, dxdt,\,\, \mathbf I_2=\int_{R^4}\rho^{3/2}
    P_\varepsilon\,\text{div}\boldsymbol\xi dxdt,\\
    \Phi=\int_{R^4}(\mathbf F\cdot\boldsymbol\xi-\partial_t\lambda\rho\mathbf v\cdot
    \boldsymbol \xi-\mathbb
    T:\nabla\boldsymbol\xi)\, dx-\int_{R^4}
    G\cdot(\rho\mathbf v)\, dxdt.
\end{split}\end{equation}

{\it  Step2. Identities.} Our next task is to calculate $\mathbf
I_1$ and $\mathbf I_2$. Note that the entries of symmetric matrix
$\nabla \boldsymbol\xi$ equal $K_{ij}*\rho$. From this and symmetry
of the kernel $K_j$ we conclude that
\begin{equation}\label{slava52}
\mathbf
I_1=\frac{1}{2}\int_{R^4}\rho(x,t)\Big\{\int_{B(x,r)}\rho(y,t)\boldsymbol
\Upsilon(x,y,t)\, dy\Big\}\,dxdt,
\end{equation}
where
$$
\boldsymbol \Upsilon(x,y,t)=K_{ij}(x-y)
\big(v_i(x,t)v_j(x,t)+v_i(y,t)v_j(y,t)-2v_i(x,t)v_j(y,t)\,\big).
$$
It is easy to see that for all  $a,b\in R^3$,
\begin{equation*}\begin{split}
\frac{|x|^2\delta_{ij}-x_ix_j}{|x|^{2}}a_ib_i=\Pi^\bot(x)
a\cdot b=\Pi^\bot(x) a\cdot\Pi^\bot(x) b, \\
\frac{x_ix_j}{|x|^2}a_jb_j=\Pi^\parallel(x) a\cdot
b=\Pi^\parallel(x) a\cdot\Pi^\parallel(x) b,
\end{split}\end{equation*}
where projections $\Pi^\bot$, $\Pi^\parallel$ are defined by
\eqref{slava1}.  From this and expression \eqref{slava51} for
$K_{ij}$ we obtain that
\begin{equation}\label{slava53}\begin{split}
\Upsilon(x,y,t)=\frac{1}{|x-y|^\mu}\big(\Pi^\bot(x-y)(\mathbf
v(x,t)-\mathbf
v(y,t))\big)^2+\\\frac{1-\mu}{|x-y|^\mu}\big(\Pi^\parallel(x-y)(\mathbf
v(x,t)-\mathbf v(y,t))\big)^2-\frac{1}{r^\mu}(\mathbf v(x,t)-\mathbf
v(y,t))^2.
\end{split}\end{equation}
for  $y\in B(x,r)$. It is clear that $\Upsilon(x,y,t)=0$ for
$|y-x|>r$.

Let us consider the term $\mathbf I_2$. We have
\begin{equation*}
    \mathbf I_2=\int_{R^4} \rho^{3/2}(x,t)P_\varepsilon(x,t)
    \text{~div~}\boldsymbol\xi\, dxdt,\,
    \text{~div~}\boldsymbol\xi=\sum_i\int_{B(x,r)}K_{ii}(x-y)\rho(y,t).
\end{equation*}
On the other hand, expression  \eqref{slava51} yields
$$
\sum_i
K_{ii}(x)=\frac{2}{|x|^\mu}+\frac{1-\mu}{|x|^\mu}-\frac{3}{r^\mu}.
$$
It follows that
\begin{equation*}\begin{split}
    \mathbf I_2=2\int_{R^4}\rho^{3/2}(x,t) P_\varepsilon(x,t)\Big\{\int_{B(x,r)}\frac
    {\rho(y,t)}{|x-y|^\mu}\, dy\big\}dxdt+\\
(1-\mu)\int_{R^4}\rho^{3/2}(x,t)
P_\varepsilon(x,t)\Big\{\int_{B(x,r)}\frac
    {\rho(y,t)}{|x-y|^\mu}\, dy\Big\}dxdt-\\\frac{3}{r^\mu}\int_{R^4}\rho^{3/2}(x,t)
P_\varepsilon(x,t)\Big\{\int_{B(x,r)}\rho(y,t)\, dy\Big\}dxdt.
\end{split}\end{equation*}
Combining this equality with \eqref{slava52},\eqref{slava53} and
recalling expressions \eqref{mila1} for the "energies" $\mathbf
E^\top$, $\mathbf E^\parallel$ we finally obtain
\begin{equation}\label{slava54}\begin{split}
    \mathbf I_1+\mathbf I_2=\int_{R^4}\rho(x,t)\Big\{\int_{B(x,r)}\Big(\frac
    {\mathbf E^\bot(x,y,t)+(1-\mu)\mathbf E^\parallel(x,y,t)}{|x-y|^\mu}
    \Big)\rho(y,t)dy\Big\}dxdt\\-
    \int_{R^4}\rho(x,t)\Big\{\int_{B(x,r}\frac{\mathbf E(x,y,t)}{r^\mu}\rho(y,t)
    \, dy\Big\}dxdt.
\end{split}\end{equation}

{\it Step 3. Estimates of $\Phi$.} Now consider the low order term
$\Phi$. Its  estimate is given by the following lemma.
\begin{lemma}\label{slava55}
Under the above assumptions,
\begin{equation}\label{slava56}
    |\Phi(r)|\leq c\sqrt{r}
\end{equation}
where $c$ is independent of $\mu$.
\end{lemma}

\begin{proof}
We begin with the observation that
\begin{multline}\label{slava56a}
\int_{R^4} (\mathbf F\cdot\boldsymbol\xi-\mathbb
T:\nabla\boldsymbol\xi)\, dxdt= \int_{R^4}\mathbf
F(x,t)\Big\{\int_{B(x,r)}
K_0(x-y)\rho(y,t)dy\Big\}dxdt\\-\int_{R^4}\rho(x,t)
\Big\{\int_{B(x,r)} K_{ij}(x-y) T_{ij}(y,t)dy\Big\}dxdt.
\end{multline}
Recall that $|K_0(x)|\leq c$ and $|K_j(x)|\leq c/|x|$. It follows
that
\begin{multline*}\begin{split}
\int_{B(x,r)} |K_0(x-y)\rho(y,t)|dy\leq\int_{B(x,r)}\rho dy\leq\\
r\big(\int_{B(x,r)}\rho^{3/2}dy\big)^{2/3}\leq
cr\|\rho\|_{L^\infty(R; L^{3/2}(R^3))}\leq cr
\end{split}\end{multline*}
and
\begin{multline*}\begin{split}
\int_{B(x,r)}
|K_{ij}(x-y)T_{ij}(y,t)|dy\leq\int_{B(x,r)}\frac{|\mathbb T|}{|x-y|}
dy\leq\\ c\sqrt{r}\big(\int_{B(x,r)} |\mathbb T|^2dy\big)^{1/2}\leq
c\sqrt{r}\|\mathbb T(t)\|_{ L^{2}(R^3))}.
\end{split}\end{multline*}
Substituting these relations in \eqref{slava56a} we finally obtain
\begin{equation}\label{slava57}\begin{split}
\Big|\int_{R^4} (\mathbf F\cdot\boldsymbol\zeta-\mathbb
T:\nabla\boldsymbol\zeta)\, dxdt\Big|\leq\\ c r\|\mathbf
F\|_{L^1(R^4)}+c\sqrt{r}\|\rho\|_{L^\infty(R; L^1(\mathbb
R))}\|\mathbb T\|_{L^2(R^4)}\leq c\sqrt{r},
\end{split}\end{equation}
Next, recall the expression for $G$
$$
G=[\lambda, K_j]*(\rho u_{\varepsilon j})+\lambda K_0*(\lambda^4_t
\varrho_\varepsilon)+\varrho_\varepsilon\mathbf u_\varepsilon\cdot
\nabla\lambda^4).
$$
Here the commutator  is defined by the equality
$$
[\lambda, K_j]*(\rho u_{\varepsilon
j})=\int_{R^3}(\lambda(x,t)-\lambda(y,t))K_j(x-y)\rho(y,t)
u_{\varepsilon j}(y,t)\, dy
$$
It is easy to see that for fixed $x$, the kernel of this operator is
compactly supported in $B(x,r)$ and is uniformly bounded. Since
$K_0(x-y)$ is uniformly bounded and is compactly supported in
$B(x,r)$, we  have the estimate
\begin{equation*}\begin{split}
|G(x,t)|\leq \int_{B(x,r)} (|\rho\mathbf u_\varepsilon|+\rho)\,
dy\leq\\
\,\Big(\int_{B(x,r)}
\varrho_\varepsilon\,dy\Big)^{1/2}\Big(\int_{B(x,r)}
\varrho_\varepsilon |\mathbf
u_\varepsilon|^2\,dy\Big)^{1/2}+\int_{B(x,r)}
\varrho_\varepsilon\,dy
\end{split}\end{equation*}
Recall that the energy estimates implies
$$
\|\varrho_\varepsilon|\mathbf
u_\varepsilon|^2\|_{L^\infty(R;L^1(R^3))}\leq c, \,\,
\int_{B(x,r)}\rho\, dy\leq
c{r}\Big(\int_{B(x,r)}\rho^{3/2}dy\Big)^{2/3}\leq cr.
$$
From this we conclude that $|G|\leq c\sqrt{r}$. It follows that
$$
\Big|\int_{R^4}
    G\cdot(\rho\mathbf v)\, dxdt\Big|\leq\int_{R^4}|\rho\mathbf v||G|\, dxdt\leq
c\sqrt{r}\int_{R^4}\rho|\mathbf v|\, dxdt\leq c\sqrt{r}.
$$
Combining this inequality  with \eqref{slava57} we obtain desired
estimate \eqref{slava56}. This completes the proof of Lemma
\ref{slava55}.
\end{proof}

\paragraph{Step 4. Finalization of the proof.}
The quantities $\mathbf I_i$ and $\Phi$ depend on the parameters
$\mu$ and $r$. Now our goal is to pass to the limit in equalities
\eqref{slava42a} and  \eqref{slava54} as $\mu\to 1$. Throughout  of
this section we write $\mathbf I_{i\mu}(r)$ and $\Phi_\mu(r)$
instead of $\mathbf I_i$ and $\Phi$. Introduce the quantities
\begin{equation}\label{slava58}\begin{split}
    \mathbf Q_\mu(r)=\int_{R^4}\rho(x,t)\Big\{\int_{B(x,r)}\big(
    \frac{\mathbf E^\bot(x,y,t)}{|x-y|^\mu}-\frac{\mathbf
    E(x,y,t)}{r^\mu}\big)\, dy\Big\}dxdt, \\ \mathbf R_\mu(r)(\mu)=
\int_{R^4}\rho(x,t)\Big\{\int_{B(x,r)}
    \frac{\mathbf E^\parallel(x,y,t)}{|x-y|^\mu}\, dy\Big\}dxdt.
\end{split}\end{equation}
With this notation equality \eqref{slava42a} becomes
\begin{equation}\label{slava59}
    \mathbf Q_\mu(r)+(1-\mu)\mathbf R_\mu(r)+\Phi_\mu(r)=0.
\end{equation}
Our  task is to  pass to the limit as $\mu\to 1$. Obviously we have
\begin{equation}\label{slava60}\begin{split}
    \Phi_\mu(r)-\frac{1}{r^\mu}\int_{R^4}\rho(x,t)\Big\{\int_{B(x,r)}
\mathbf
    E(x,y,t)\, dy\Big\}dxdt\to\\
\Phi_1(r)-\frac{1}{r}\int_{R^4}\rho(x,t)\Big\{\int_{B(x,r)} \mathbf
    E(x,y,t)\, dy\Big\}dxdt\text{~as~}\mu\to 1.
\end{split}\end{equation}
It follows from expression \eqref{mila1} for $\mathbf E^\parallel$
that
$$
\mathbf E^\parallel\leq |\mathbf v(x,t)|^2+|\mathbf
v(y,t)|^2+\rho(x,t)^{1/2} P_\varepsilon
$$
Thus we get
\begin{equation}\label{slava61}\begin{split}
    \mathbf R_\mu(r)\leq\\
    \int_{R^4}\rho(x,t)\Big\{\int_{B(x,r)}\big( |\mathbf
    v(y,t)|^2+|\mathbf
    v(x,t)|^2+\rho(x,t)^{1/2}P_\varepsilon(x,t)\big)\frac{\rho(y,t)dy}{|x-y|^\mu}
    \Big\}dxdt=\\
 \int_{R^4}\big(\rho(x,t)|\mathbf
    v(x,t)|^2+\rho(x,t)^{3/2}P_\varepsilon(x,t)\big)
    \Big\{\int_{B(x,r)}\frac{\rho(y,t)dy}{|x-y|^\mu}
    \Big\}dxdt+\\\int_{R^4}\rho(x,t)
    \Big\{\int_{B(x,r)}\varrho(y,t)|\mathbf v(y,t)|^2\frac{dy}{|x-y|^\mu}
    \Big\}dxdt=\\
\int_{R^4}\big(2\rho(x,t)|\mathbf
    v(x,t)|^2+\rho(x,t)^{3/2}P_\varepsilon(x,t)\big)
    \Big\{\int_{B(x,r)}\frac{\rho(y,t)dy}{|x-y|^\mu}
    \Big\}dxdt.
\end{split}\end{equation}
On the other hand, we have
\begin{multline*}
\int_{B(x,r)}\frac{\rho(y,t)dy}{|x-y|^\mu}\leq
\Big(\int_{B(x,r)}\frac
{dy}{|x-y|^{3\mu}}\Big)^{1/3}\Big(\int_{B(x,r)}\varrho(y,t)^{3/2}dy
\Big)^{2/3}\leq\\ c(1-\mu)^{-1/3}\|\varrho^{3/2}\|_{L^\infty(R;
L^1(R^3))}\leq c(1-\mu)^{-1/3}.
\end{multline*}
Next, the energy estimate implies
$$
\|\rho |\mathbf v|^2+\rho^{3/2}+\rho^{3/2}
P_\varepsilon\|_{L^\infty(R; L^1(R^3))}\leq c,
$$
which along with \eqref{slava61} yields
$$
(1-\mu)\mathbf R_\mu(r)\leq c(1-\mu)^{2/3}\to 0\text{~~as~~}\mu\to
1.
$$
Recall expression \eqref{slava58} for $\mathbf Q_\mu$.  Letting
$\mu\to 1$ in equality \eqref{slava59} and using relation
\eqref{slava60} we obtain
\begin{multline*}
    \lim\limits_{\mu\to
    1}\int_{R^4}\rho(x,t)\Big\{\int_{B(x,r)}\mathbf
    E^\bot(x,y,t)\frac{\rho(y,t)\,
    dy}{|x-y|^\mu}\Big\}dxdt\\-\frac{1}{r}\int_{R^4}\Big\{\int_{B(x,r)}
    \mathbf E(x,y,t)\, dy\Big\}dxdt+\Phi_1(r)=0.
\end{multline*}
Note that for $r\leq 1$, the integrand in the first integral is a
monotone function of $\mu$. It follows from this and the Fatou
theorem that
\begin{equation*}\begin{split}
\lim\limits_{\mu\to
    1}\int_{R^4}\rho(x,t)\Big\{\int_{B(x,r)}\mathbf
    E^\bot(x,y,t)\frac{\rho(y,t)\,
    dy}{|x-y|^\mu}\Big\}dxdt=\\
     \int_{R^4}\rho(x,t)\Big\{\int_{B(x,r)}\mathbf
    E^\bot(x,y,t)\frac{\rho(y,t)\,
    dy}{|x-y|}\Big\}dxdt.
\end{split}\end{equation*}
Thus we get
\begin{multline}\label{slava62}
    \int_{R^4}\rho(x,t)\Big\{\int_{B(x,r)}\mathbf
    E^\bot(x,y,t)\frac{\rho(y,t)\,
    dy}{|x-y|}\Big\}dxdt-\\\frac{1}{r}\int_{R^4}\Big\{\int_{B(x,r)}
   \mathbf E(x,y,t)\, dy\Big\}dxdt +\Phi(r),
\end{multline}
with $\Phi=\Phi_1$. Let us transform the first integral in the left
hand side of this equality. It is easy to check that
\begin{multline*}
\int_{B(x,r)}\mathbf
    E^\bot(x,y,t)\frac{\rho(y,t)\,
    dy}{|x-y|}= \int_0^r\Big\{\int_{|x-y|=s}
    \mathbf E^\bot(x,y,t)\rho(y,t)\,dS\Big\}\frac{ds}{s}=   \\
\int_0^r\Big\{\int_{B(x,s)}
    \mathbf E^\bot(x,y,t)\rho(y,t)\,dy\Big\}\frac{ds}{s^2}+\frac{1}{r}\int_{B(x,r)}
    \mathbf E^\bot(x,y,t)\rho(y,t)\,dy.
\end{multline*}
Substituting this equality into \eqref{slava62} and recalling
expressions \eqref{slava9} and \eqref{slava10} for $\Gamma$ and $W$
we arrive at desired identity \eqref{slava30}. It remains to prove
estimates \eqref{slava14} and \eqref{slava15}. Estimate
\eqref{slava15} obviously follows from formula  estimate
\eqref{slava56} in Lemma \ref{slava55}.   In order to estimate
$W(r)$, notice that the energy estimates \eqref{mambda3},
\eqref{mambda4} and expression \eqref{slava9} for $\Gamma$ imply the
inequality $\Gamma(1)\leq c$. From this, identity \eqref{slava30},
and estimate \eqref{slava15} for $\Phi$ we obtain
$$
W(r)\leq W(1)\leq \Gamma(1)+|\Phi(1)|\leq c\text{~~for~~} r\in
(0,1],
$$
which yields desired estimate \eqref{slava14} for $W$. It remains to
note that estimate \eqref{slava14} for $\Gamma$ obviously  follows
from estimate \eqref{slava15}, identity \eqref{slava30}, and
estimate for $W$.

%%%%%%%%%%%%%%%%%%%%%%%%%%%%%%%%%%%%%%%%%%%%%%%%%%%%%%%%%%%%%%%%%%%%%%%%%%%%%%%%%%%%%%%%%
%%%%%%%%%%%%%%%%%%%%%%%%%%%%%%%%%%%%%%%%%%%%%%%%%%%%%%%%%%%%%%%%%%%%%%%%%%%%%%%%%%%%%%%%%

\section{ Estimates  of $W$ and $\Gamma$. }\label{sasha}
\renewcommand{\theequation}{4.\arabic{equation}}
\setcounter{equation}{0} In this section we derive auxiliary
estimates for the functions $\Gamma$ and $W$ defined by equalities
\eqref{slava9} and \eqref{slava10}. The following proposition is the
main result of this section. Recall definition \eqref{slava10} of
the function $ W(r)$
\begin{proposition}\label{sasha05}
Whenever $\mu<4/3$,
\begin{equation}\label{sasha06}
\Gamma(r)+W(r)\leq c|\log
    r|^{-\mu}\text{~~for all~~} r\in (0, e^{-1}).
\end{equation}
\end{proposition}
The rest of the section is devoted to the proof of  Proposition
\ref{sasha05}. Our strategy is based on the analysis of identity
\eqref{slava30} in Proposition \ref{slava14p}. Rewrite this identity
in the form
\begin{equation}\label{sasha07}
    \frac{r}{a(r)}\, W'(r)-W(r)=\Phi(r),
    \text{~~where~~}a(r)=\frac{rW'(r)}{\Gamma(r)}.
\end{equation}
We will consider  \eqref{sasha07} as a differential equation on the
interval $(0, e^{-1})$. If the coefficient $a$ oscillates as $r\to
0$, then the behavior of solution may be complicated. In order to
cope with this difficulty, we split the interval $(0, e^{-1}]$ into
two parts corresponding "large" and "small" values of the
coefficient $a(r)$. Fix an arbitrary $C>5$ and introduce the sets
\begin{equation}\label{sasha1}\begin{split}
  A=\big\{ r\in (0, e^{-1}]:\,\, a(r)> C/ |\ln r|\,\big\},\\
B=\big\{ r\in (0, e^{-1}]:\,\, a(r)\leq C/ |\ln r|\,\big\}.
\end{split}\end{equation}
It is clear that $ (0, e^{-1}]=\, A\, \cup \, B. $ Next set
\begin{equation}\label{sasha2}
\begin{split}
    W_A(r)\,=\,\frac{1}{2}\int_{[0,r]\cap A}
    \int_{ R^4}\rho(x,t)\frac{1}{s^2}\int_{B(x,s)}
    (\Pi^\bot(x-y)\,(\mathbf v(x,t)-\mathbf v(y,t)))^2\,dy\,dx\,dt\,ds\\
    +\,2\int_{[0,r]\cap A}\int_{ R^4}\rho^{\frac{3}{2}}(x,t)P_\varepsilon\,\frac{1}{s^2}
    \int_{B(x,s)}\rho(y,t)\,dy\,dx\,dt\,ds,
\end{split}
\end{equation}
and
\begin{equation}\label{sasha3}
W_B(r)=W(r)-W_A(r)
\end{equation}
It is clear that $W_B$ is defined by equality \eqref{sasha2} with
$A$ replaced by $B$. Now set
\begin{equation}\label{sasha4}
    \Gamma_A=\chi_A(r)\Gamma(r), \quad \Gamma_B=\chi_B(r)\Gamma(r),
\end{equation}
where $\chi_A$ and $\chi_B$ are characteristic functions of $A$ and
$B$.

\subsection{ Estimate of  $W_B$. }\label{victor}

In this section we obtain the   estimate for the function $W_B(r)$
given by formula \eqref{slava10}. Our goal is to prove that
$W(r)\sim |\ln r|^{-\mu}$ for all sufficiently small $r$ for every
positive $\mu\in [0, 4/3)$.
 The result is given by the following
\begin{lemma}\label{victor1} For every $\mu\in [0,4/3)$, there
is a constant $c>0$, depending only on $C$ and $\mu$, such that
\begin{equation}\label{victor2}
 \Gamma_B(\sigma)+   W_B(\sigma)\leq \frac{c}{|\ln \sigma|^{\mu}}
 \text{~~for all ~~} \sigma \in (0,
    e^{-1}).
\end{equation}

\end{lemma}
\begin{proof}
We begin with the observation that definition  \eqref{sasha1} of the
set $B$ and formula  \eqref{sasha07} for the coefficient $a$ imply
$$
rW'(r)\leq \frac{C}{|\log r|} \Gamma(r)\text{~~for~~} r\in B.
$$
Next note that  definition \eqref{slava10} for $W$ and  formulae
\eqref{mila1}, \eqref{slava9} we obtain
\begin{equation}\begin{split}\label{victor47}
\frac{1}{2}\int_{ R^4}\rho(x,t)\frac{1}{r}
    \int_{B(x,r)}\rho(y,t)\,(\Pi^\bot(x-y)\,(\mathbf v(x,t)-\mathbf v(y,t)))^2
    \,dy\,dx\,dt
    \\
    +\,2\int_{ R^4}\rho^{\frac{3}{2}}(x,t)P_\varepsilon\,
    \frac{1}{r}\int_{B(x,r)}\rho(y,t)\,dy\,dx\,dt\\
    =\int_{R^4}\rho(x,t)\int_{B(x,r)} \rho(y,t)\mathbf E^\bot(x,y,t)dy
    =rW'(r).
\end{split}\end{equation}
Recalling formulae \eqref{slava9} for $\Gamma$ and  \eqref{mila1}
for $\mathbf E^\parallel$ we arrive at the inequality
\begin{equation*}\begin{split}
\frac{1}{2}\int_{R^4}\rho(x,t)\frac{1}{r}
    \int_{B(x,r)}\rho(y,t)\,(\Pi^\bot(x-y)\,(\mathbf v(x,t)-\mathbf v(y,t)))^2
    \,dy\,dx\,dt
    \\
    +\,2\int_{ R^4}\rho^{\frac{3}{2}}(x,t)P_\varepsilon\,
    \frac{1}{r}\int_{B(x,r)}\rho(y,t)\,dy\,dx\,dt\leq \\\frac{C}{|\ln r|}\Big\{
\frac{1}{2}\int_{
R^4}\rho(x,t)\,\frac{1}{r}\int_{B(x,r)}\rho(y,t)\,(\Pi^\parallel(x-y)\,
    (\mathbf v(x,t)-\mathbf v(y,t)))^2\,dy\,dx\,dt\,\\
    +\int_{ R^4}\rho^{\frac{3}{2}}(x,t)P_\varepsilon\,
    \frac{1}{r}\int_{B(x,r)}\rho(y,t)\,dy\,dx\,dt\Big\}
\end{split}\end{equation*}
which holds true for $r\in (0, e^{-1}]\cap B$. Now choose $\tau_0$
so small that $C/|\log r|\leq 1$ for all $r\in (0,\tau_0)$. From
this and \eqref{victor47} we conclude  that the inequality
\begin{equation}\begin{split}\label{victor41}
r W'(r)\leq 2 \frac{C}{|\ln r|}\Big\{ \frac{1}{2}\int_{
R^4}\rho(x,t)\,\frac{1}{r}\int_{B(x,r)}\rho(y,t)\,(\Pi^\parallel(x-y)\,
    (\mathbf v(x,t)-\mathbf v(y,t)))^2\,dy\,dx\,dt\,
\end{split}\end{equation}
holds for all $r\in B\cap (0,\tau_0)$. On the other hand, we have
the inequality
\begin{equation*}\begin{split}
\int_{R^4}\rho(x,t)\,\frac{1}{r}\int_{B(x,r)}\rho(y,t)\,(\Pi^\parallel(x-y)\,
    (\mathbf v(x,t)-\mathbf v(y,t)))^2\,dy\,dx\,dt\leq\\
    c\int_{
R^4}\rho(x,t)\,\frac{1}{r}\int_{B(x,r)}\rho(y,t)|\mathbf
v(y,t)^2\,\,dy\,dx\,dt +c\int_{ R^4}\rho(x,t)\mathbf
v(x,t)^2\,\frac{1}{r}\int_{B(x,r)}\rho(y,t)\,\,dy\,dx\,dt\\=2c\int_{
R^4}\rho(x,t)\mathbf
v(x,t)^2\,\frac{1}{r}\int_{B(x,r)}\rho(y,t)\,\,dy\,dx\,dt.
\end{split}
\end{equation*}
which being substituted in \eqref{victor41} gives
\begin{equation}\begin{split}\label{victor42}
rW'(r)\leq  \frac{c}{|\ln r|}\int_{ R^4}\rho(x,t)|\mathbf
v(x,t)|^2\,\frac{1}{r}\int_{B(x,r)}\rho(y,t)\,\,dy\,dx\,dt
\end{split}\end{equation}
for all $r\in B\cap (0,\tau_0)$.  Now fix an arbitrary $N>0$ and set
\begin{equation}\label{victor43}
\mathcal E_N=\big\{\, t\in \mathbb R:\,\,\, \|\mathbf
v(\cdot,t)\|_{W^{1,2}(\mathbb R^3)}\leq N\,\big\}.
\end{equation}
It follows from the Young inequality that
\begin{gather*}
\frac{c}{|\log r|}\int_{ R^4}\rho(x,t)|\mathbf
v(x,t)|^2\,\frac{1}{r}\int_{B(x,r)}\rho(y,t)\,\,dy\,dx\,dt \leq\\
 \int_{\mathcal E_N}\int_{
R^3}\rho(x,t)^{3/2}\,\frac{1}{r}\int_{B(x,r)}\rho(y,t)\,\,dy\,dx\,dt
+\frac{c}{|\log r|^3}\int_{\mathcal E_N}\int_{ R^3}|\mathbf
v (x,t)|^6\,\frac{1}{r}\int_{B(x,r)}\rho(y,t)\,\,dy\,dx\,dt+\\
\frac{c}{|\log r|}\int_{R\setminus \mathcal E_N}\int_{
R^3}\rho(x,t)\mathbf
v(x,t)^2\,\frac{1}{r}\int_{B(x,r)}\rho(y,t)\,\,dy\,dx\,dt
\end{gather*}
Substituting this inequality in the right hand side of
\eqref{victor42},  noting that $P_\varepsilon\geq 1$, and recalling
formula \eqref{victor47} for $rW'(r)$ we obtain for all $r\in
B\cap(0, \tau_0]$
\begin{equation}\begin{split}\label{victor44}
rW'(r)\leq \frac{c}{|\log r|^3}\int_{\mathcal E_N}\int_{
R^3}|\mathbf
v (x,t)|^6\,\frac{1}{r}\int_{B(x,r)}\rho(y,t)\,\,dy\,dx\,dt+\\
\frac{c}{|\log r|}\int_{R\setminus \mathcal E_N}\int_{
R^3}\rho(x,t)\mathbf
v(x,t)^2\,\frac{1}{r}\int_{B(x,r)}\rho(y,t)\,\,dy\,dx\,dt
\end{split}\end{equation}
Now our task is to estimate the integrals in the right hand side of
this inequality. We begin with the observation that for every
$\alpha\in (0,1)$, we have
$$
\int_{B(x,
e^{-1})}\frac{\rho(y,t)dy}{|x-y||\log|x-y||^\alpha}=\int_0^{e^{-1}}\Big\{\int_{
B(x,r)}\rho(y,t)\, dy\Big\}\frac{dr}{r|\log r|^\alpha}.
$$
Next note that
\begin{multline*}
\int_0^{e^{-1}}\Big\{\int_{\partial B(x,r)}\rho(y,t)\,
dy\Big\}\frac{dr}{r|\log
r|^\alpha}=\int_0^{e^{-1}}\Big\{\int_{B(x,r)}\rho(y,t)\,
dy\Big\}(1+\frac{\alpha}{|\log r|})\frac{dr}{r^2|\log r|^\alpha}\\
+e\int_{\partial B(x,e^{-1})}\rho(y,t)\, dy.
\end{multline*}
Since $|\log r|\geq 1$ for $r\in (0, e^{-1})$, it follows from this
that
\begin{equation*}
\int_0^{e^{-1}}\Big\{\int_{B(x,r)}\rho(y,t)\,
dy\Big\}\frac{dr}{r^2|\log r|^\alpha}\leq c\int_{B(x,
e^{-1})}\frac{\rho(y,t)dy}{|x-y||\log|x-y||^\alpha}.
\end{equation*}
Now fix $\alpha\in (1/3, 1)$. We have
\begin{multline*}
\int_{B(x, e^{-1})}\frac{\rho(y,t)dy}{|x-y||\log|x-y||^\alpha}
\leq\\
\Big(\int_{B(x,
e^{-1})}\frac{dy}{|x-y|^3|\log|x-y||^{3\alpha}}\Big)^{1/3}
\Big(\int_{B(x, e^{-1})}\rho^{3/2}(y,t)dy\Big)^{2/3}\leq c(\alpha)
\end{multline*}
Combining the obtained result we conclude that the inequality
\begin{equation}\label{victor46}
\int_0^{e^{-1}}\varkappa(x,t,r)\frac{dr}{r}\leq c(\alpha), \quad
\varkappa(x,t,r)= \frac{1}{r|\log r|^\alpha}\int_{B(x,r)}\rho(y,t)\,
dy
\end{equation}
holds for all $\alpha\in (1/3, 1)$. Here the constant $c$ depends
only on $\alpha$. Substituting this identity into \eqref{victor44}
and using estimate \eqref{victor46} we conclude that the inequality
\begin{equation}\begin{split}\label{victor48}
r W'(r)\leq  \frac{c}{|\log r|^{3-\alpha}}\int_{\mathcal E_N}\int_{
R^3}\varkappa(x,t,r)|\mathbf
v (x,t)|^6\,dx\,dt+\\
\frac{c}{|\log r|^{1-\alpha}}\int_{R\setminus \mathcal E_N}\int_{
R^3}\varkappa(x,t,r)\rho(x,t)\mathbf v(x,t)^2\,dx\,dt.
\end{split}\end{equation}
holds true for all $r\in B\cap (0, \tau_0]$.

Notice that
\begin{equation}\begin{split}\label{victor209}
\int_{\mathcal E_N}\int_{R^3} |\mathbf v|^6\,dxdt\leq
c\int_{\mathcal E_N}\|\mathbf v(\cdot,t)\|_{W^{1,2}(R^3)}^6\,dt\leq\\
c N^4\int_{\mathcal E_N}\|\mathbf
v(\cdot,t)\|_{W^{1,2}(R^3)}^2\,dt\leq c N^4
\end{split}\end{equation}
On the other hand, the energy estimate implies
\begin{equation}\begin{split}\label{victor210}
\int_{R\setminus\mathcal E_N}\int_{R^3} \rho|\mathbf v|^2\,dxdt\leq
c\text{~~meas~}([0,T]\setminus \mathcal E_N)\leq c N^{-2},
\end{split}\end{equation}
since $\rho|\mathbf v|^2$ vanishes for $t\in R\setminus [0,T]$.

Now choose an arbitrary $\sigma\in (0, e^{-1})$. It follows from
that for every $r\in B\cap (0,\sigma]$, we have
\begin{equation*}\begin{split}
 W'(r)\leq  \frac{c}{|\log \sigma|^{3-\alpha}}\int_{\mathcal E_N}\int_{
R^3}\frac{\varkappa(x,t,r)}{r}|\mathbf
v (x,t)|^6\,dx\,dt+\\
\frac{c}{|\log \sigma|^{1-\alpha}}\int_{R\setminus \mathcal
E_N}\int_{ R^3}\frac{\varkappa(x,t,r)}{r}\rho(x,t)\mathbf
v(x,t)^2\,dx\,dt.
\end{split}\end{equation*}
Integrating both sides of this inequality with respect to $r$ over
the set $B\cap (0,\sigma]$ and using estimates \eqref{victor46},
\eqref{victor209}, and \eqref{victor210} we arrive at the inequality
\begin{equation*}
W_B(\sigma)=\int_{B\cap [0,\sigma]}W'(r)\,dr\leq c\frac{N^4}{|\log
\sigma|^{3-\alpha}}+\frac{c}{N^2|\log \sigma|^{1-\alpha}}.
\end{equation*}
Choosing $N^6=|\log\sigma|^2$ we finally obtain
$$
W_B{\sigma}\leq c |\log\sigma|^{-5/3+\alpha}
$$
where $\alpha\in (1/3,1)$ is an arbitrary number.  This gives the
estimate \eqref{victor2} for $W_B$ with $\mu=5/3-\alpha$.
%%%%%%%%%%%%%%%%%%%%%%%%%%%%%%%%%%%%%%%%%%%%%%%%%%%%%%%%%%%%%%%%%%%%%%%%%%%%%%%%%%%%%%%%%
%%%%%%%%%%%%%%%%%%%%%%%%%%%%%%%%%%%%%%%%%%%%%%%%%%%%%%%%%%%%%%%%%%%%%%%%%%%%%%%%%%%%%%%%%
%%%%%%%%%%%%%%%%%%%%%%%%%%%%%%%%%%%%%%%%%%%%%%%%%%%%%%%%%%%%%%%%%%%%%%%%%%%%%%%%%%%%%%%%%%%%
%%%%%%%%%%%%%%%%%%%%%%%%%%%%%%%%%%%%%%%%%%%%%%%%%%%%%%%%%%%%%%%%%%%%%%%%%%%%%%%%%%%%%%%%%%%%%

The proof of estimate for
 $\Gamma_B$ is similar. Recalling expression  \eqref{slava9} for $\Gamma(r)$ and
 arguing as in the proof of  \eqref{victor44} we arrive at the
 inequality
\begin{equation}\begin{split}\label{victor44a}
\Gamma(r)\leq \frac{c}{|\log r|^3}\int_{\mathcal E_N}\int_{
R^3}|\mathbf
v (x,t)|^6\,\frac{1}{r}\int_{B(x,r)}\rho(y,t)\,\,dy\,dx\,dt+\\
\frac{c}{|\log r|}\int_{R\setminus \mathcal E_N}\int_{
R^3}\rho(x,t)\mathbf
v(x,t)^2\,\frac{1}{r}\int_{B(x,r)}\rho(y,t)\,\,dy\,dx\,dt
\end{split}\end{equation}
which holds true for all $r\in B\cap (0,\tau_0)$.   Here the set
$\mathcal E_N$ is defined by the equality

Now our task is to estimate the integrals in the right hand side of
this inequality. We begin with the observation that for every
$\alpha\in (0,1)$ and $r\in (0, e^{-1}$, we have
$$
\int_{B(x,r)}\frac{\rho(y,t)dy}{|x-y||\log|x-y||^\alpha}=\int_0^{r}\Big\{\int_{\partial
B(x,s)}\rho(y,t)\, dy\Big\}\frac{ds}{s|\log s|^\alpha}.
$$
Next note that
\begin{multline*}
\int_0^{r}\Big\{\int_{\partial B(x,s)}\rho(y,t)\,
dy\Big\}\frac{ds}{s|\log
s|^\alpha}=\int_0^{r}\Big\{\int_{B(x,s)}\rho(y,t)\,
dy\Big\}(1+\frac{\alpha}{|\log s|})\frac{ds}{s^2|\log s|^\alpha}\\
+\frac{1}{r|\log r|^\alpha}\int_{\partial B(x,r)}\rho(y,t)\, dy.
\end{multline*}
It follows from this that
\begin{equation*}
\frac{1}{r|\log r|^\alpha}\int_{\partial B(x,r)}\rho(y,t)\, dy\leq
c\int_{B(x, e^{-1})}\frac{\rho(y,t)dy}{|x-y||\log|x-y||^\alpha}.
\end{equation*}
Now fix $\alpha\in (1/3, 1)$. We have
\begin{multline*}
\int_{B(x, e^{-1})}\frac{\rho(y,t)dy}{|x-y||\log|x-y||^\alpha}
\leq\\
\Big(\int_{B(x,
e^{-1})}\frac{dy}{|x-y|^3|\log|x-y||^{3\alpha}}\Big)^{1/3}
\Big(\int_{B(x, e^{-1})}\rho^{3/2}(y,t)dy\Big)^{2/3}\leq c(\alpha)
\end{multline*}
Combining the obtained result we conclude that the inequality
\begin{equation}\label{victor46a}
\frac{1}{r|\log r|^\alpha}\int_{\partial B(x,r)}\rho(y,t)\, dy\leq
c(\alpha)
\end{equation}
holds for all $\alpha\in (1/3, 1)$ and $r\in (0, e^{-1})$. Here the
constant $c$ depends only on $\alpha$. Substituting this inequality
into \eqref{victor44a} we obtain that the inequality
\begin{equation}\begin{split}\label{victor48a}
\Gamma(r)\leq  \frac{c}{|\log r|^{3-\alpha}}\int_{\mathcal
E_N}\int_{ R^3}|\mathbf
v (x,t)|^6\,dx\,dt+\\
\frac{c}{|\log r|^{1-\alpha}}\int_{R\setminus \mathcal E_N}\int_{
R^3}\rho(x,t)\mathbf v(x,t)^2\,dx\,dt,
\end{split}\end{equation}
holds true for all $r\in B\cap (0, \tau_0]$. Notice that
\begin{equation*}\begin{split}
\int_{\mathcal E_N}\int_{R^3} |\mathbf v|^6\,dxdt\leq
c\int_{\mathcal E_N}\|\mathbf v(\cdot,t)\|_{W^{1,2}(R^3)}^6\,dt\leq\\
c N^4\int_{\mathcal E_N}\|\mathbf
v(\cdot,t)\|_{W^{1,2}(R^3)}^2\,dt\leq c N^4
\end{split}\end{equation*}
On the other hand, the energy estimate implies
\begin{equation*}\begin{split}
\int_{R\setminus\mathcal E_N}\int_{R^3} \rho|\mathbf v|^2\,dxdt\leq
c\text{~~meas~}([0,T]\setminus \mathcal E_N)\leq c N^{-2},
\end{split}\end{equation*}
since $\rho|\mathbf v|^2$ vanishes for $t\in R\setminus [0,T]$. From
this and \eqref{victor48a} we conclude that the inequality
\begin{equation*}
W_B(\sigma)=\int_{B\cap [0,\sigma]}W'(r)\,dr\leq c\frac{N^4}{|\log
\sigma|^{3-\alpha}}+\frac{c}{N^2|\log \sigma|^{1-\alpha}}.
\end{equation*}
holds for all $r\in B\cap (0,\tau_0)$. Choosing $N^6=|\log r|^2$ we
finally obtain
$$
\Gamma{r}\leq c |\log r|^{-5/3+\alpha}\text{~~for~~} r\in B\cap
(0,\tau_0),
$$
where $\alpha\in (1/3,1)$ is an arbitrary number. It remains  to
note that $\Gamma_B(r)= \Gamma(r)$ for $r\in B$ and $\Gamma_B(r)=0$
otherwise. This gives the estimate \eqref{victor2} for $\Gamma_B$
with $\mu=5/3-\alpha$. This completes the proof of Lemma
\ref{victor1}
\end{proof}
\subsection{Estimate of $W_A$ and $\Gamma_A$.
Proof of Proposition \ref{sasha05}}\label{sashaA} Our next task is
to estimate $W_A$. The result is given by the following lemma
\begin{lemma}\label{victor51} For every $\mu\in [0,4/3)$, there
is a constant $c>0$, depending only on $C$ and $\mu$, such that
\begin{equation}\label{victor52}
 \Gamma_A(\sigma)+   W_A(\sigma)\leq \frac{c}{|\ln \sigma|^{\mu}}
 \text{~~for all ~~} \sigma \in (0,
    e^{-1}).
\end{equation}

\end{lemma}
\begin{proof}
Recall that
\begin{equation*}
    \Gamma_A=\chi_A \Gamma=\chi_A
    \frac{r}{a}W'(r)=\frac{r}{a}W'_A(r)\text{~~and hence~~}
    \Gamma=\frac{r}{a}W'_A(r)+\Gamma_B.
\end{equation*}
 Substituting these equalities into
\eqref{sasha07} we arrive at the equality
\begin{equation}\label{sasha6}
  \frac{r}{a}W'_A(r)-W_A(r) =-\Gamma_B(r)+W_B(r)+\Phi(r), \quad r\in (0, e^{-1}).
\end{equation}
Since
$$
a(r)\geq \frac{C}{|\log r|}\text{~~on~~} A\text{~~and~~}
W'_A(r)=0\text{~~on~~} (0,e^{-1})\setminus A,
$$
it follows that
\begin{equation}\label{sasha7}
  \frac{C}{|\log r|}  W'_A(r)-W_A(r) \geq-\Gamma_B(r)-|\Phi(r)|.
\end{equation}
It follows from estimate \eqref{victor2} in Proposition
\ref{victor1} and estimate \eqref{slava15} in Proposition
\ref{slava14p} that
\begin{equation}\label{victor19}
    |\Gamma_B(r)|\leq c|\log r|^{-2/3}, \quad |\Phi(r)|\leq c \sqrt{r}
    \text{~~for ~~} r\in (0, e^{-1}),
\end{equation}
Since $C>4$, the function  $W_A$ satisfies all conditions of Lemma
\ref{Alemma} with $\beta_i=\mu$, $R=\tau=e^{-1}$. Applying this
lemma  we obtain desired estimate
\begin{equation}\label{victor20}
    W_A(r)\leq c |\log r|^{-\mu} \text{~~for~~} r\in (0, e^{-1}).
\end{equation}
for the function $W_A$. The estimate for $\Gamma_A$ follows from the
estimate for $W_B$, inequality \eqref{victor2}, estimate
\eqref{slava15} for $\Phi$, and the identity
$$
\Gamma_A=-\Gamma_B +W_B+W_A+\Phi.
$$
This completes the proof of Lemma \ref{victor51}.
\end{proof}

Finally notice that Proposition \ref{sasha05} is the straightforward
consequence of Lemmas \ref{victor1} and \ref{victor51}.
\subsection{ Preliminary
potential estimates} In this section we employ Proposition
\ref{sasha05} in order to estimate a potential of the density
function. The result is given by
\begin{proposition}\label{victor53}
For every $\mu\in (0,4/3)$ there is $c>0$ such that
\begin{equation}\label{frank16}
    \int_{\mathbb R^4}\varrho^{3/2}(x,t)
    \Big\{\int_{\mathbb R^3}\frac{\big|\log |x-y|\big|^{\mu}+1}
    {|x-y|}\varrho(y,t)dy\Big\}dxdt\leq c.
\end{equation}
\end{proposition}
\begin{proof} Since the function $\rho$ is compactly supported in
$Q_T$ and $\|\rho\|_{L^\infty(R; L^{3/2}(R^3))}\leq c$, it suffices
to prove that
\begin{equation}\label{frank30}
\int_{\mathbb R^4}\varrho^{3/2}(x,t)
    \Big\{\int_{B(x, e^{-1})}\frac{\big|\log |x-y|\big|^{\mu}}
    {|x-y|}\varrho(y,t)dy\Big\}dxdt\leq c.
\end{equation}
It is easy to see that
\begin{equation*}\begin{split}
  \int_{B(x, e^{-1})}\frac{\big|\log |x-y|\big|^\mu}{|x-y|}
   \varrho(y,t)dy =\int_0^{e^{-1}}\frac{|\log r|^\mu}{r}
  \Big\{\int_{\partial B(x,r)}\rho(y,t)dS\Big\}dr
  =\\
  \int_0^{e^{-1}} \big(1+\frac{\mu}{|\log r|}\big)\frac{|\log r|^\mu}{r^2}
\Big\{\int_{ B(x,r)}\rho(y,t)dy\Big\}dr+e\int_{
 B(x,e^{-1})}\rho(y,t)dy\leq\\
 c\int_0^{e^{-1}} \frac{|\log r|^\mu}{r^2}
\Big\{\int_{ B(x,r)}\rho(y,t)dy\Big\}dr+c,
\end{split}\end{equation*}
which along with formula \eqref{slava10} yields
\begin{equation}\label{frank101}\begin{split}
\int_0^{e^{-1}}\int_{R^4}\rho^{3/2}(x,t)\Big\{\int_{B(x,
e^{-1})}\frac{\big|\log |x-y|\big|^\mu}{|x-y|}
   \varrho(y,t)dy\Big\}\, dxdt\leq  \\
c\int_0^{e^{-1}}\int_{R^4}\rho^{3/2}(x,t)\Big\{\int_{
B(x,r)}\rho(y,t)dy\Big\}\frac{|\log r|^\mu}{r^2}\, dxdtdr+c\leq\\
c\int_0^{e^{-1}}|\log r|^\mu\, W'(r)dr+c\leq
c\int_0^{e^{-1}}\frac{|\log r|^{\mu-1}}{r}\, W(r)dr+c.
\end{split}\end{equation}
Now set
$$
\alpha=\frac{1}{2}\Big(\,\mu+\frac{4}{3}\,\Big), \quad \beta=1+
\frac{1}{2}\Big(\,\frac{4}{3}-\mu\,\Big), \quad \alpha\in (0,
\frac{4}{3}), \quad \beta>1, \quad \alpha-\beta=\mu-1.
$$
With this notation we can rewrite inequality \eqref{frank101} in the
equivalent form
\begin{equation}\label{frank102}\begin{split}
\int_0^{e^{-1}}\int_{R^4}\rho^{3/2}(x,t)\Big\{\int_{B(x,
e^{-1})}\frac{\big|\log |x-y|\big|^\mu}{|x-y|}
   \varrho(y,t)dy\Big\}\, dxdt\leq  \\
 c\int_0^{e^{-1}}\frac{1}{r|\log r|^\beta}\, \big(|\log r|^\alpha W(r)\big)dr+c.
\end{split}\end{equation}
By virtue of Proposition \ref{sasha05}, we have $|\log r|^\alpha
W(r)\leq c$. From this and \eqref{frank102} we conclude that
\begin{equation*}\begin{split}
\int_0^{e^{-1}}\int_{R^4}\rho^{3/2}(x,t)\Big\{\int_{B(x,
e^{-1})}\frac{\big|\log |x-y|\big|^\mu}{|x-y|}
 \varrho(y,t)dy\Big\}\, dxdt\leq  \\
c\int_0^{e^{-1}}\frac{1}{r|\log r|^\beta}\,dr+c\leq c,
\end{split}\end{equation*}
and the proposition follows.

\end{proof}

\section{Nonlinear potentials}\label{frank}
\renewcommand{\theequation}{5.\arabic{equation}}
\setcounter{equation}{0} In this section we estimate various
integral operators of the density function.  In particular we will
consider the following potentials
\begin{gather}\label{frank3}
 \Psi_\mu(x,t)=\int_{
 R^3}\frac{1}{|x-y|^2}(1+|\ln|x-y||)^{\mu}\varrho(y,t)^{4/5}dy,\\
 \label{flint6}
    G_\mu(x)=\int_{ R^3}\frac{1}{|x-y|^2}
    \varrho^{5/4}(y,t)\ln (1+\varrho(y,t))^\mu\, dy,\\
\label{flint7}
    \Sigma_\mu(x)=G_0 (1+|\ln G_0|)^\mu, \quad G_0=
    \int_{ R^3}\frac{1}{|x-y|^2} \varrho^{5/4}(y,t)\, dy.
\end{gather}
The goal of this section is the proof of the following proposition.
\begin{proposition} \label{frank1} For every $\mu\in (0,2/3)$, there
is a constant $c>0$ such that
\begin{gather}\label{frank2}
    \| \Psi_\mu\|_{L^2(R^4)}\leq c\\
\label{flint4}
    \|G_\mu\|_{L^2( R^4)}\leq c,\\
\label{flint5}
    \|\Sigma_\mu\|_{L^2( R^4)}\leq c.
\end{gather}
\end{proposition}
The proof falls into three steps.

\paragraph{Step 1.}First we prove
estimate \eqref{frank2}. Our consideration is   based on the
following elementary lemma.
\begin{lemma}\label{frank17}
For every $\mu\geq 0$,
\begin{equation}\label{frank14}\begin{split}
I:=\int_{ R^3}\frac{|\ln |x-z||^{\mu/2}}{|x-z|^2}\, \frac{|\ln
|y-z||^{\mu/2}}{|y-z|^2}\,dz\leq c \frac{|\ln |x-y||^{\mu}}{|x-y|}
\end{split}\end{equation}
\end{lemma}
\begin{proof} The proof is given in Section \ref{petr}. \end{proof}

Now fix an arbitrary $\mu\in (0,2/3) $ and set $\lambda=2\mu<4/3$.
We have
\begin{gather*}
\int_{ R^7}\frac{|\ln |x-y||^\lambda+1}{|x-y|}
\varrho^{5/4}(x,t)\varrho^{5/4}(y,t) dxdydt\leq \\\int_{ R^7}
\frac{|\ln |x-y||^\lambda+1}{|x-y|} \varrho(x,t)\varrho(y,t)
(\sqrt{\varrho}(x,t)+
\sqrt{\varrho}(x,t))dxdydt=\\
\int_{ R^7}\frac{|\ln |x-y||^\lambda+1}{|x-y|}
(\varrho^{3/2}(x,t)\varrho(y,t) +
\varrho(x,t)\varrho^{3/2}(y,t))dxdydt=\\
2\int_{ R^7} \varrho^{3/2}\Big\{\int_{\mathbb R^3}\frac{|\ln
|x-y||^\lambda+1}{|x-y|} \varrho(y,t)dy\Big\}dxdt.
\end{gather*}
 From this   and Proposition \ref{victor53} we conclude  that
\begin{equation}\label{frank20}
    \int_{ R^7}\frac{|\ln |x-y||^\lambda+1}{|x-y|}
     \varrho^{5/4}(x,t)\varrho^{5/4}(y,t) dxdydt\leq c.
\end{equation}
Now we employ  Lemma \ref{frank17} to obtain
\begin{equation*}\begin{split}
\int_{ R^4}\Big\{\int_{\mathbb R^3} \frac{|\ln
|x-z||^{\lambda/2}}{|x-z|^2} \varrho^{5/4}(x,t)dx\Big\}
\Big\{\int_{\mathbb R^3}\frac{|\ln |y-z||^{\lambda/2}}{|y-z|^2} \varrho^{5/4}(y,t)dy\Big\}dz=\\
\int_{ R^7}\Big\{\int_{\mathbb R^3}\frac{|\ln
|x-z||^{\lambda/2}}{|x-z|^2} \frac{|\ln
|y-z||^{\lambda/2}}{|y-z|^2}dz\Big\}
\varrho^{5/4}(x,t)\varrho^{5/4}(y,t) dxdydt\leq \\c \int_{
R^7}\frac{|\ln |x-y||^\lambda+1}{|x-y|}
\varrho^{5/4}(x,t)\varrho^{5/4}(y,t) dxdydt\leq c.
\end{split}\end{equation*}
Since $\lambda=2\mu$, we arrive at the inequality
\begin{equation}\begin{split}\label{frank21}
\int_{ R^4}\Big\{\int_{\mathbb R^3}\frac{|\ln |x-y||^{\mu}}{|x-y|^2}
\varrho^{5/4}(y,t)dy\Big\}^2dxdt\leq c.
\end{split}\end{equation}
Repeating these arguments we arrive at the estimate
\begin{equation*}
\int_{\mathbb R^3\times \mathbb R}\Big\{\int_{\mathbb
R^3}\frac{1}{|x-y|^2} \varrho^{5/4}(y,t)dy\Big\}^2dxdt\leq c.
\end{equation*}
Combining this result with estimate \eqref{frank21} we obtain
desired inequality \eqref{frank2}.
 \paragraph{Step 2.} Let us turn to the proof of estimate
 \eqref{flint4}. Fix an arbitrary $\sigma\in (\mu, 2/3)$.
 It follows from the  Cauchy
 inequality that
\begin{gather*}
G_\mu\leq \Big(\int_{ R^3}\frac{|\ln |x-y||^\sigma+1}{|x-y|^2}
\varrho^{5/4}(y,t)\, dy\Big)^{1/2} \Big(\int_{\mathbb R^3}\frac{\ln
(1+\rho)^{2\mu}}{|x-y|^2(|\ln |x-y||^\sigma+1)} \varrho^{5/4}(y,t)\,
dy\Big)^{1/2}\\\leq \Psi_\sigma+\int_{\mathbb R^3}\frac{\ln
(1+\rho)^{2\mu}}{|x-y|^2(|\ln |x-y||^\sigma+1)} \varrho^{5/4}(y,t)\,
dy\leq\Psi_\sigma+ \\c\Big(\int_{ R^3}\frac{|\ln
|x-y||^\sigma+1}{|x-y|^2} \varrho^{5/4}(y,t)\, dy\Big)^{1/2}
\Big(\int_{ R^3}\frac{\ln (1+\rho)^{4\mu}}{|x-y|^2(|\ln
|x-y||^\sigma+1)^3} \varrho^{5/4}(y,t)\, dy\Big)^{1/2}\\\leq
2\Psi_\sigma+\int_{ R^3}\frac{\ln (1+\rho)^{4\mu}}{|x-y|^2(|\ln
|x-y||^\sigma+1)^3}
 \varrho^{5/4}(y,t)\, dy\leq\\
2\Psi_\sigma+c\int_{ R^3}\frac{\ln (1+\rho)^{4\mu}}{|x-y|^2(|\ln
|x-y||+1)^{3\sigma }} \varrho^{5/4}(y,t)\, dy.
\end{gather*}
Repeating these arguments we conclude that for any integer $n>1$,
\begin{equation}\label{flint8}
    G_\mu\leq n\Psi_\sigma+c(n) I_n, \text{~~where~~} I_n=\int_{\mathbb R^3}
    \frac{\ln (1+\rho)^{\mu_n}}{|x-y|^2(|\ln |x-y||+1)^{\sigma_n}} \varrho^{5/4}(y,t)\, dy,
\end{equation}
where
\begin{equation}\label{flint9}
    \mu_n=2^n\mu, \quad \sigma_n=(2^n-1)\sigma, \quad \sigma_{n+1}=
    2\sigma_n+1,\quad n\geq 2.
\end{equation}
Now our task is to estimate $I_n$ for all large $n$. Further we will
assume that $n$ is sufficiently large number, which will be
specified below. Choose an arbitrary $\varepsilon \in (0,1)$.
Applying the H\'{o}lder inequality we obtain
\begin{equation*}\begin{split}
 I_n \leq \Big(\int_{ R^3}
 \frac{dy}{|x-y|^3(|\ln|x-y||+1)^{1+\varepsilon}}
 \varrho^{3/2}(y,t)\, dy\Big)^{1/2}
 \Big(\int_{ R^3}\varrho^{3/2}(y,t)dy\Big)^{\frac{5}{6}
 -\frac{1}{2}}\times\\
 \times \Big(\int_{ R^3}
 \frac{dy}{|x-y|^3(|\ln|x-y||+1)^{6(\sigma_n-\frac{1+\varepsilon}{2})}}
 |\ln (1+\varrho)|^{6\mu_n}\Big)^{1/6}.
\end{split}\end{equation*}
From this and the energy estimate we conclude that
\begin{equation}\begin{split}\label{flint10}
 I_n \leq c\Big(\int_{ R^3}
 \frac{1}{|x-y|^3(|\ln|x-y||+1)^{1+\varepsilon}}
 \varrho^{3/2}(y,t)\, dy\Big)^{1/2}
 \times\\
 \times \Big(\int_{ R^3}
 \frac{dy}{|x-y|^3(|\ln|x-y||+1)^{6(\sigma_n-\frac{1+\varepsilon}{2})}}
 |\ln (1+\varrho)|^{6\mu_n}\Big)^{1/6}.
\end{split}\end{equation}
Let us estimate the second  integral in the right hand side. By
abuse of notation we set
$$
k_n=6(\sigma_n-\frac{1+\varepsilon}{2})\to\infty, \quad
\beta_n=6\mu_n\to \infty \text{~~as~~} n\to\infty,
$$
and
\begin{equation}\begin{split}\label{flint11}
 J_n=\int_{ R^3}\frac{1}{|x-y|^3(|\ln|x-y||+1)^{k_n}}
 |\ln (1+\varrho)|^{\beta_n}dy.
\end{split}\end{equation}
Introduce the Young function $ \Psi(z)=|z| \ln (1+|z|)^{2k_n/3}$.
Since $k_n\geq 2$, it is convex. Obviously the conjugate function
$\Psi^*$ admits the estimate
$$
0\leq \Phi^*(z)\leq c (e^{|z|^\frac{3}{2k_n}}-1).
$$
The Young inequality yields the estimate
\begin{equation*}\begin{split}
 J_n\leq
 \int_{ R^3} \frac{1}{|x-y|^3(|\ln|x-y||+1)^{k_n}}
 \big(\ln (1+\frac{1}{|x-y|^3(|\ln|x-y||+1)^{k_n}}\big)^{2k_n/3}\, dy+\\
 c\int_{ R^3}
 \exp\big(\ln(1+\varrho)^{\frac{3}{2k_n}})-1\big)\, dy.
 \end{split}\end{equation*}
Next notice that
\begin{equation*}\begin{split}
\ln (1+\frac{1}{|x-y|^3(|\ln|x-y||+1)^{k_n}}\big)\leq \ln
(1+\frac{1}{|x-y|^3})\\
 +\ln(1+\frac{1}{(|\ln|x-y||+1)^{k_n}})\leq c (1+|\ln|x-y||),
 \end{split}\end{equation*}
which gives
\begin{equation}\begin{split}\label{flint11a}
 J_n\leq
 \int_{ R^3} \frac{1}{|x-y|^3(|\ln|x-y||+1)^{k_n/3}}dy+\\
 c\int_{ R^3}
 \big(\exp\big(\ln(1+\varrho)^{\frac{3\beta_n}{2k_n}}\big)-1\big)\, dy
 \end{split}\end{equation}
Note that
$$\beta_n/k_n\to \mu/ \sigma=q<1,\quad k_n\to\infty\text{~~as~~} n\to \infty.
$$
Hence we may fix $n$ such that
$$
0<\beta_n/k_n<1,\quad k_n>3
$$
For such $n$, we have
$$
0\leq \exp\big(\ln(1+\varrho)^{\frac{3\beta_n}{2k_n}}\big)-1\leq c
\varrho^{3/2}+c \varrho.
$$
It follows that
\begin{equation*}\begin{split}
\int_{ R^3}
\big(\exp\big(\ln(1+\varrho)^{\frac{\beta_n}{k_n}}\big)-1\big)\,
dy\leq c,\\
   \int_{\mathbb R^3} \frac{1}{|x-y|^3(|\ln|x-y||+1)^{k_n/3}}dy\leq c.
\end{split}\end{equation*} Combining this result with
\eqref{flint11a} we finally obtain  $J_n\leq c$. Substituting this
estimate into \eqref{flint10} we arrive at the inequality
\begin{equation*}\begin{split}
 I_n \leq c\Big(\int_{ R^3}\frac{1}{|x-y|^3(|\ln|x-y||+1)^{1+\varepsilon}}\varrho^{3/2}(y,t)\, dy\Big)^{1/2}.
\end{split}\end{equation*}
From this we obtain
\begin{equation*}\begin{split}
    \int_{ R^3} I_n^2(x,t)\,dx\leq c \int_{\mathbb R^6}
    \frac{1}{|x-y|^3(|\ln|x-y||+1)^{1+\varepsilon}}\varrho^{3/2}(y,t)\, dydx=\\
 \int_{ R^3}\Big\{\int_{\mathbb R^3}
 \frac{1}{|x-y|^3(|\ln|x-y||+1)^{1+\varepsilon}}dx\Big\}\varrho^{3/2}(y,t)\, dy
 \leq
 c\int_{\mathbb R^3}\varrho^{3/2}(y,t)\, dy\leq c.
\end{split}\end{equation*}
On the other hand, formula \eqref{flint8} for $I_n$ implies that
$I_n(x,t)$ vanishes for $t\in R\setminus  [0,T]$. From this,
identity \eqref{flint8}, and estimate \eqref{frank2} we obtain
desired inequality \eqref{flint4}.

\paragraph{Step 3.} It remains to prove estimate \eqref{flint5}.
Let us consider a convex function $g(z)=|z|\ln(1+|z|)^\mu.$   Since
$\rho$ is supported in $\Omega\times [0,T]$, it follows from the
Jensen inequality  that
\begin{equation}\begin{split}\label{flint14}
G_0\ln (1+G_0)^\mu= g(G_0)\leq c\int_\Omega g\Big(
\frac{1}{|x-y|^2}\rho^{5/4}(y,t)\Big)\, dy.
\end{split}\end{equation}
Note that
\begin{equation*}\begin{split}
 g\Big( \frac{1}{|x-y|^2} \varrho^{5/4}(y,t)\, \Big)& \leq  \frac{1}{|x-y|^2} \varrho^{5/4}
 \ln (1+ \frac{1}{|x-y|^2})^\mu\\&+\frac{1}{|x-y|^2} \varrho^{5/4}\ln
 (1+\varrho^{5/4})^\mu,
\end{split} \end{equation*} which along with the obvious inequalities
$$
 \ln (1+ \frac{1}{|x-y|^2})^\mu\leq c(1+ |\ln|x-y||)^\mu,\quad
\ln (1+\varrho^{5/4})^\mu\leq c \ln(1+\varrho)^\mu,
$$
implies
\begin{equation*}\begin{split}
     g\Big( \frac{1}{|x-y|^2} \varrho^{5/4}(y,t)\, \Big) \leq c
      \frac{1+|\ln|x-y||^\mu}{|x-y|^2}\varrho^{5/4}+
      \frac{1}{|x-y|^2} \varrho^{5/4}\ln(1+\varrho)^\mu.
\end{split}\end{equation*}
Substituting this estimate  into  \eqref{flint14}and recalling
formulae \eqref{frank3}, \eqref{flint6},  we obtain
\begin{equation*}
    \int_{\Omega}
g\Big( \frac{1}{|x-y|^2} \varrho^{5/4}(y,t)\, \Big)dy\leq c\Psi_\mu
+ cG_\mu
\end{equation*}
It remains to notice that now  desired estimate \eqref{flint5} is a
consequence of estimates \eqref{frank2} and \eqref{flint4}.

\section{Proof of Theorem  \ref{pont7}}\label{zirda}
\renewcommand{\theequation}{6.\arabic{equation}}
\setcounter{equation}{0}

The proof falls into two steps.  First we derive an auxiliary
estimate for the total kinetic energy  as a function of the temporal
variable.
\paragraph{Auxiliary estimates.}
Fix an arbitrary $\mu\in (0,2/3)$ and set $\alpha=4\mu/5$. Introduce
the functions
\begin{equation}\label{zirda1}
\mathbf H (t)=    \int_{ R^3}\rho(t)|\mathbf v(t)|^2\,
    \log(1+\rho(t)|\mathbf v(t)|^2)^\alpha\, dx, \quad
    \mathbf N(t)=\|\mathbf
    v(\cdot, t)\|_{W^{1,2}(R^3)}.
\end{equation}
Our first task is to estimate $\mathbf H(t)$ in terms of $\mathbf
N(t)$. The corresponding result is given by the following
\begin{proposition}\label{zirda1p}
Under the above assumptions there is a nonnegative function $F(t)$
such that
\begin{equation}\label{zirda14}
\int_{R}F(t)\, dt\, \leq c, \quad \mathbf H(t)\leq F(t)\,\mathbf
N(t)^{5/4}.
\end{equation}
\end{proposition}

\begin{proof}
By abuse of notation further we omit the dependence on $t$. It
follows from the H\"{o}lder inequality that
\begin{equation}\label{zirda2}
    \mathbf H\leq \Big(\int_{ R^3} \rho^{5/4} |\mathbf
    v|\log(1+\rho |\mathbf v|^2)^{\mu}\, dx\Big)^{4/5}
\Big(\int_{ R^3}|\mathbf v|^6\, dx\Big)^{1/5}.
\end{equation}
Next set
\begin{equation}\label{zirda3}
\mathbf v=\mathbf N \mathbf w, \quad \|\mathbf w\|_{W^{1,2}(
R^3)}=1.
\end{equation}
We have
\begin{gather*}
\log(1+\rho |\mathbf v|^2)^{\mu}\leq c
\log(1+\rho)^{\mu}+c\log(1+ |\mathbf v|)^{\mu},\\
 \log(1+ |\mathbf v|^2)^{\mu}\leq c\log(1+\mathbf
N)^{\mu}+c\log(1+|\mathbf w|)^{\mu},\\
\Big(\int_{ R^3}|\mathbf v|^6\, dx\Big)^{1/5}=\mathbf N^{6/5}
\Big(\int_{ R^3}|\mathbf w|^6\, dx\Big)^{1/5}\leq c \mathbf N^{6/5}.
\end{gather*}
Substituting these inequalities into \eqref{zirda2} we arrive at the
inequality
\begin{equation}\label{zirda6}
\mathbf H\leq \mathbf N^{6/5}( \mathbf H_1^{4/5}+\mathbf H_2^{4/5}+
\mathbf H_3^{4/5}),
\end{equation}
where
\begin{equation}\label{zirda4}\begin{split}
\mathbf H_1=&\mathbf N\int_{ R^3} \rho^{5/4}|\mathbf
w|\log(1+|\mathbf w|)^{\mu}\,
dx,\\
\mathbf H_2=& \mathbf N \log(1+\mathbf N)^{\mu}\int_{ R^3}
\rho^{5/4}|\mathbf w|\,
dx,\\
\mathbf H_3=&\mathbf N\,\int_{ R^3} \rho^{5/4} \log(1+\rho)^{\mu}
|\mathbf w|\, dx.
\end{split}\end{equation}
Now our task is to estimate $\mathbf H_k$. The result is given by
the following lemma
\begin{lemma}\label{zirdalemma}
Under the assumptions of Proposition \ref{zirda1p}, we have
\begin{equation}\label{zirda125}
    \|\mathbf H_i\|_{L^1(R)}\leq c, \quad i=1,2,3.
\end{equation}
\end{lemma}
\begin{proof}
The proof falls into three steps.\\
\noindent {\it Step 1.} We start with the estimating of $\mathbf
H_1$. Introduce the vector field
$$
\boldsymbol \zeta=\nabla \big(|\mathbf w|\log(1+|\mathbf
w|)^\mu\big).
$$
It follows from the integral representation of Sobolev functions
that
\begin{equation}\label{zirda50}
\mathbf H_1\leq c\mathbf N\int_{ R^3}\Big\{ \int_{
R^3}\frac{\rho(y,t)^{5/4}}{|x-y|^2}\, dy\Big\}\,
|\boldsymbol\zeta(x,t)|\,dx=c\int_{ R^3}G_0(x,t)\,
|\boldsymbol\zeta(x,t)|\,dx.
\end{equation}
where $G_0$ is given by \eqref{flint7}. In order to estimate
$\boldsymbol \zeta$, note that that the function
$g(z)=|z|\log(1+|z|)^\mu$ is convex and its congjugate admits the
estimate
$$
g^*(z)\leq c(\exp\{z^{1/\mu}\}-1)
$$
From this and the Young inequality we obtain the estimate
\begin{equation*}\begin{split}
|\boldsymbol\zeta|&\leq c (\log(1+|\mathbf w|)^\mu +1)|\nabla
\mathbf w|\\&\leq  c |\nabla \mathbf w|\log(1+| \nabla \mathbf
w|)^\mu+ c\exp\{(\log(1+|\mathbf w|)^\mu+1)^{1/\mu}\}-1.
\end{split}\end{equation*}
Combining this result with the inequality
$$
\exp\{(\log(1+|\mathbf w|)^\mu+1)^{1/\mu}\}-1\leq c |\mathbf w|+1
$$
and noting that for every $t$, the vector field $\boldsymbol
\zeta(\cdot, t)$ is compactly supported in $\Omega$ we arrive  at
the estimate
$$
|\boldsymbol \zeta|\leq c|\nabla \mathbf w|\log(1+|\nabla \mathbf
w|)^\mu + c |\mathbf w|+c\boldsymbol{1}_\Omega=c g(|\nabla \mathbf
w|)+c| \mathbf w|+c\boldsymbol{1}_\Omega.
$$
Thus we get
\begin{equation}\label{zirda7}\begin{split}
\mathbf H_1\leq c\mathbf N\int_{ R^3} G_0\,\,g(|\nabla \mathbf w|)
\,dx+c\mathbf N\int_{ R^3} G_0\,\big(| \mathbf
w|+\boldsymbol{1}_\Omega\,\big)\, dx.
\end{split}\end{equation}
Notice that the function $\Upsilon(z)=z^2\log (e+|z|)^{2\mu}$ is
convex and its conjugate admits the estimate
\begin{equation}\label{zirda30}
\Upsilon^*(z)\leq c |z|^2\log(e+|z|)^{-2\mu}
\end{equation}
The Young inequality implies that
\begin{equation}\label{zirda30a}
G_0 \,\,g\leq \Upsilon\big(\,\mathbf N^{-1/2}G_0\,\big)+
\Upsilon^*\big(\,\mathbf N^{1/2}g\,\big)
\end{equation}
For $\mathbf N\geq 1$ we have
\begin{equation*}\begin{split}
 \Upsilon\big(\,\mathbf N^{-1/2}G_0\,\big)\leq  \mathbf N^{-1}
 G_0^2\log(e+\mathbf N^{-1/2}G_0)^{2\mu}\leq\\ c \mathbf N^{-1}G_0^2\log
 (1+G_0)^{2\mu}+c \mathbf N^{-1}G_0^2=\mathbf N^{-1} \Sigma_\mu^2 +
 +c \mathbf N^{-1}G_0^2,
 \end{split}\end{equation*}
where $\Sigma_\mu$ is defined by equality \eqref{flint7}. On the
other hand, we have
\begin{equation*}\begin{split}
\Upsilon^*\big(\,\mathbf N^{1/2} g(z)\,\big)\leq \mathbf N
g^2\log(e+\mathbf N^{1/2}
g)^{-2\mu}\leq\\
\mathbf N z^2\log(e+ |z|)^{2\mu} \log\big(e+\mathbf
N^{1/2}|z|\log(e+|z|)^\mu\big)^{-2\mu}\leq c\mathbf N z^2.
\end{split}\end{equation*}
Substituting estimates for $\Upsilon$ and $\Upsilon^*$ in
\eqref{zirda30a} we arrive at the inequality
\begin{equation}\label{zirda21}
G_0\,\, g(\,|\nabla\mathbf w|\,)\leq c \mathbf N^{-1}( \Sigma_\mu^2
 +G_0^2)+c\mathbf N |\nabla\mathbf w|^2.
\end{equation}
Recall that  $g(\,|\nabla\mathbf w|\,)$ is compactly supported in
$Q_T=\Omega\times (0,T)$. It follows from this and \eqref{zirda21}
 that for any $\mathbf N$,
\begin{equation*}
G_0\,\, g,|\nabla\mathbf w|\,)\leq c (1+\mathbf N)^{-1}(
\Sigma_\mu^2
 +G_0^2)+c\boldsymbol{1}_{Q_T}\,(1+\mathbf N) |\nabla\mathbf w|^2
\end{equation*}
 Combining this inequality with \eqref{zirda7} we arrive at the
estimate
\begin{equation}\label{zirda22a}
\mathbf H_1\leq c\int_{R^3}( \Sigma_\mu^2 +
 +G_0^2)\, dx+\boldsymbol{1}_{[0,T]}\int_\Omega \Big((1+\mathbf
 N)^2|\nabla \mathbf w|^2+\mathbf N G_0(|\mathbf w|+1)\Big)\,dx
\end{equation}
Since $\mathbf w(t)$ is supported in $\Omega$ and $\|\mathbf
w\|_{W^{1,2}(R^3)}\leq 1$, we have
$$
\int_\Omega |\nabla\mathbf w|^2\, dx\leq 1, \quad \int_\Omega
 G_0\,\,(|\mathbf w|+1)\,dx \leq \|G_0\|_{L^2(\Omega)}.
 $$
Substituting these inequalities into \eqref{zirda22a} we obtain the
estimate
\begin{equation*}
\mathbf H_1\leq c\int_{R^3}( \Sigma_\mu^2
 +G_0^2)\, dx+c\boldsymbol{1}_{[0,T]}(1+\mathbf
 N)^2+\mathbf N \|G_0\|_{L^2(\Omega)}.
\end{equation*}
which along with the Cauchy inequality implies
\begin{equation*}
\mathbf H_1\leq c\int_{R^3}( \Sigma_\mu(t)^2
 +G_0(t)^2)\, dx+c\boldsymbol{1}_{[0,T]}(1+\mathbf
 N)^2.
\end{equation*}
From this, inequalities  \eqref{flint4}, \eqref{flint5} in Lemma
\ref{frank1}, and the estimate
$$\|\mathbf N\|_{L^2(R)}\leq \|\mathbf v \|_{L^2(R;
W^{1,2}(R^3))}\leq c
$$
we conclude that
\begin{equation}\label{zirda22}
\int_{R}\mathbf H_1(t)\, dt\leq c.
\end{equation}

{\it Step 2.} Now our task is to estimate
 $\mathbf H_2$.
Notice that
\begin{equation*}\begin{split}
\int_{\mathbb R^3}\rho^{5/4} |\mathbf w|\, dx\leq c\int_{ R^3}\Big\{
\int_{\mathbb R^3}\frac{\rho(y)^{5/4}}{|x-y|^2}\, dy\Big\}\, |\nabla
\mathbf w|\,dx\leq\\ c\int_{\Omega} G_0|\nabla \mathbf w|\, dx\leq
\Big(\int_{\Omega}G_0^2\, dx\Big)^{1/2},
\end{split}\end{equation*}
which gives
\begin{equation}\label{zirda10}
\mathbf H_2\leq \mathbf N\log(1+\mathbf N)^{\mu}\sqrt{\psi},
\text{~~where~~} \psi=\int_{\Omega}G_0^2\,dx.
\end{equation}
The Young inequality implies
$$
\mathbf N\log(1+\mathbf N)^{\mu}\sqrt{\psi}\leq c
\Upsilon\big(\,\sqrt{\psi}\,\big)+c \Upsilon^*\big(\,\mathbf
N\log(1+\mathbf N)^{\mu}\,\big),
$$
where $\Upsilon(z)=z^2\log(e+|z|)^{2\mu}$. From this and
\eqref{zirda30} we obtain
\begin{equation}\label{zirda25}\begin{split}
\mathbf H_2\leq c\mathbf N\log(1+\mathbf N)^{\mu}\sqrt{\psi}\leq c
\frac{\mathbf N^2\log(1+\mathbf N)^{2\mu}}{\log\big(\,e+ \mathbf
N\log(1+
\mathbf N)^{\mu}\,\big)^{2\mu}}+ \psi\log(e+\sqrt{\psi})^{2\mu}\\
\leq c \mathbf N^2+c\psi\log(e+\psi)^{2\mu}\leq c\mathbf
N^2+c\psi\log(1+\psi)^{2\mu}+c\psi.
\end{split}\end{equation}
Since the function $|z|\log(1+|z|)^{2\mu}$ is convex, we have
$$
\psi\log(1+\psi)^{2\mu}\leq c \int_\Omega
G_0^2\log(1+G_0^2)^{2\mu}\, dx \leq c \int_\Omega
G_0^2\log(1+G_0)^{2\mu}\, dx=c\int_\Omega \Sigma_\mu^{2}\, dx
$$
 Substituting this inequality in \eqref{zirda25} and integrating
 the result with respect to  $t$ we finally obtain
\begin{equation}\label{zirda26}\begin{split}
\int_R \mathbf H_2(t)\, dt \leq c\int_R \mathbf N(t)^2\,
dt+\int_{R^4}\Sigma_\mu ^2\, dxdt+\int_{R^4}G_0^2\, dxdt\leq c.
\end{split}\end{equation}
{\it Step 3.} It remains to estimate $\mathbf H_3$. Since $\|\mathbf
w(t)\|_{W^{1,2}(R^3)}=1$, we have
\begin{equation*}\begin{split}
|\mathbf H_3|\leq \mathbf N c\int_{\mathbb R^3}\Big\{ \int_{\mathbb
R^3}\frac{\rho(y)^{5/4}\log(1+\rho)^\mu}{|x-y|^2}\, dy\Big\}\,
|\nabla w|\,dx\leq \\c\mathbf N\int_{\mathbb R^3} G_\mu|\nabla w|\,
dx\leq c\mathbf N\Big(\int_{\mathbb R^3}G_\mu^2\, dx\Big)^{1/2}.
\end{split}\end{equation*}
From this and estimate \eqref{flint4} we conclude that
\begin{equation}\label{zirda27}\begin{split}
\int_R \mathbf H_3(t)\, dt \leq c\int_R \mathbf N(t)^2\,
dt+\int_{R^4}G_\mu ^2\, dxdt\leq c.
\end{split}\end{equation}
This completes the proof of Lemma \ref{zirdalemma}
\end{proof}
Let us turn to the proof of Proposition \ref{zirda1p}.   Set
$$
F(t)= \mathbf N(t)^{2/5}\big(\mathbf H_1(t)^{4/5}+\mathbf
H_2(t)^{4/5}+\mathbf H_3(t)^{4/5}\big).
$$
It remains to note  that the desired estimate \eqref{zirda14} for
$F$ obviously follows from \eqref{zirda22}, \eqref{zirda26},
\eqref{zirda27}, and the H\"{o}lder inequality, and the proposition
follows.
\end{proof}

\paragraph{Proof of Theorem \ref{pont7}.}
 We are now in a position to complete the proof of Theorem
 \ref{pont7}.
 For every $s>0$, denote by $E(s)$ the
measure of the set $\{\rho|\mathbf v|^2\geq s\}$, i.e.,
$$
E(s)=\text{~meas~}\big\{(x,t):\,\, \rho(x,t)|\mathbf v(x,t)|^2\geq
s\big\}
$$
The function $E$ is positive, decreasing, and tends to $0$ as $s\to
\infty$.  Introduce the   Stiltjies measure $\boldsymbol\mu(s)=
-dE(s)$. For every nonnegative, decreasing, continuously
differentiable on $(0, \infty )$, and continuous on $[0,\infty)$
function $f(s)$ with $f(0)=0$, we have
\begin{equation}\label{tantra1}
    \int_{R^4}f(\rho|\mathbf v|^2)\, dxdt=\int_{(0,\infty)}
    f(s)\, d\boldsymbol\mu(s).
\end{equation}
This equality holds true if the integral in the left hand side
exists. Now fix an arbitrary $M>0$ and set
\begin{gather*}
\mathcal T=\big\{t\in  R\,:\,\,\, \|\mathbf
v(t)\|_{W^{1,2}( R^3)}\,<\,M\,\big\},\\
 E^+(s)=\text{~meas~}\big\{(x,t):\,\, \rho(x,t)|\mathbf
v(x,t)|^2\geq s, \,\,\, t\in  R\setminus\mathcal T\big\},\\
 E^-(s)=\text{~meas~}\big\{(x,t):\,\, \rho(x,t)|\mathbf
v(x,t)|^2\geq s,\,\,\,  t\in \mathcal T\big\}
\end{gather*}
 It is clear that
\begin{equation*}
E(s)=E^-(s)+E^+(s), \quad
\boldsymbol\mu(s)=\boldsymbol\mu^-(s)+\boldsymbol\mu^+(s),\quad
\boldsymbol \mu^{\pm}=-d E^{\pm}.
\end{equation*}
and
\begin{equation}\begin{split}\label{tantra2}
\int_{\mathcal T}\int_{R^3}\rho|\mathbf v|^2\log(1+\rho|\mathbf
v|^2)^\alpha\, dxdt=\int_{(0,\infty)}
s\log(1+s)^\alpha\, d\boldsymbol \mu^-(s),\\
\int_{ R\setminus \mathcal T}\int_{ R^3}\rho|\mathbf v|^2\,
dxdt=\int_{(0,\infty)} s\, d\boldsymbol \mu^+(s),
\end{split}\end{equation}
Formula \eqref{zirda1} and estimate \eqref{zirda14} imply the
inequality
$$
\int_{\mathcal T}\int_{ R^3}\rho|\mathbf v|^2\log(1+\rho|\mathbf
v|^2)^\alpha\, dxdt\equiv \int_{\mathcal T} \mathbf H(t)\, dt\leq c
M^{4/5}\text{~~for~~} \alpha=4\mu/5\in (0, 8/15).
$$
 Since the function $\|\mathbf v(t)\|_{W^{1,2}(R^3)}$ is integrable
with square, we have $\text{~meas~}( R\setminus \mathcal T)\leq c
M^{-2}.$ From this and the  boundedness of $\rho|\mathbf v|^2$ in
the space $L^\infty( R; L^1(\mathbb R^3))$ we get
$$
\int_{ R\setminus\mathcal T}\int_{\mathbb R^3}\rho|\mathbf v|^2\,
dxdt\leq c M^{-2}.
$$
Combining the obtained estimates with \eqref{tantra2} we  arrive at
the inequalities
\begin{equation*}\begin{split}
\int_{(0,\infty)} s\log(1+s)^\alpha\, d\boldsymbol \mu^-(s)\leq c
M^{4/5}, \quad \int_{(0,\infty)} s\, d\boldsymbol \mu^+(s)\leq c
M^{-2}.
\end{split}\end{equation*}
In particular, for any $z>0$ we have
\begin{equation*}\begin{split}
\int_{(z,\infty)} s\, d\boldsymbol \mu^-(s)\leq c
M^{4/5}\log(1+z)^{-\alpha}, \quad \int_{(z,\infty)} s\, d\boldsymbol
\mu^+(s)\leq c M^{-2},
\end{split}\end{equation*}
which along with  \eqref{tantra3} leads to the inequality
\begin{equation}\begin{split}\label{tantra3}
\int_{(z,\infty)} s\, d\boldsymbol \mu=\int_{(z,\infty)} s\,
d\boldsymbol \mu^-+\int_{(z,\infty)} s\, d\boldsymbol \mu^+\leq c
M^{4/5}\log(1+z)^{-\alpha}+c{M^-2}.
\end{split}\end{equation}
Here the constant $c$ is independent of $M$ and $z$. Now fix an
arbitrary $z$ and set $M=\log(1+z)^{5\alpha/14}$. We have
$$
M^{-2}= M^{4/5}\log(1+z)^{-\alpha}=\log(1+z)^{-5\alpha/7}.
$$
In this case, inequality \eqref{tantra3} yields the estimate
\begin{equation}\begin{split}\label{tantra4a}
\int_{(z,\infty)} s\, d\boldsymbol \mu\leq c \log(1+z)^{-5\alpha/7}.
\end{split}\end{equation}
Notice that for every $0<\sigma<5\alpha/7=4\mu/7<8/21$, inequality
\eqref{tantra4a} implies
\begin{equation*}\begin{split}
\int_{(0,\infty)}z\log(1+z)^\sigma\, d\boldsymbol \mu(z)=\sigma
\int_{(0,\infty)}\frac{\log(1+z)^{\sigma-1}}{1+z}\Big\{\int_z^\infty
s\, d\boldsymbol\mu(s)\,\Big\}\,dz\leq\\
\int_{(0,\infty)}\frac{\log(1+z)^{\sigma-1-5\alpha/7}}{1+z}\, dz\leq
c
\end{split}\end{equation*}
Combining this result with  identity \eqref{tantra1} we obtain the
desired estimate \eqref{stokes5abc}. This completes the proof of
Theorems \ref{pont7}.
\begin{appendix}
\section{Appendix}\label{petr}
\renewcommand{\theequation}{A1.\arabic{equation}}
\setcounter{equation}{0}

\paragraph{Differential inequality.}In this section we derive estimates of solutions to a degenerate
differential inequality.  Fix positive constants  $C>5$  and $R\in
(0,e^{-1})$. Let us introduce functions $\Gamma, F:[0, R]\to \mathbb
R^+$ such that
\begin{equation}\label{petr1}\begin{split}
|\Gamma(r)|\,\leq\frac{k}{|\ln r|^{\mu}}, \quad \mu\in (0, C),\quad
|F|\,\leq\,c_1\sqrt{r}.
\end{split}\end{equation}
\begin{lemma}\label{Alemma}
Let $W_A:[0, R]\to \mathbb R^+$ satisfy the inequalities
\begin{equation}\label{petr5}
    \frac{r|\ln r|}{C}\,W_A'(r) - W_A(r)\,\geq\,-\Gamma(r)-F(r),
    \quad  |W_A(e^{-1})|\,\leq\,c_1.
\end{equation}
Then there is a  constant $c_2$, depending on $C$, $c_1$, $R$, such
that
\begin{equation}\label{petr7}
    W_A(r)\,\leq\,c_2\,(|\ln r|^{-C}+
    k\,|\ln r|^{-\mu})\text{~~for all~~}r\in
    (0,R].
\end{equation}
\end{lemma}
\begin{proof} Set
\begin{equation}\label{petr8}
    b(r)=\exp\Big({\int_r^{R} \frac{C}{|\ln s|s}\,ds}\Big)=|\ln r|^C
    |\log R|^{-C}
\end{equation}
Multiplying both sides of \eqref{petr5} by $Cb/(r|\log r|)$ and
integrating the result over $[r, R]$ we arrive at the inequality
\begin{equation}\label{petr9}
    b(r)[W_A(r)]-W_A(R)\,\leq\,C\int_r^{R}b(s)\,\frac{1}{s|\ln s|}\,(\Gamma(s)+F(s))\,ds,
\end{equation}
which yields the estimate
\begin{equation*}
    W_A(r)\,\leq\,b(r)^{-1}\,W_A(R) + C\,b(r)^{-1}\,\int_r^{R}b(s)
    \,\frac{1}{s|\ln s|}\,(\Gamma(s)+F(s))\,ds\text{~~for~~} r\in
    [0,R].
\end{equation*}
From this and \eqref{petr1}, \eqref{petr5}, \eqref{petr8}  we
conclude that
\begin{equation*}\begin{split}
    W_A(r)\leq\frac{c_1|\log R|^C}{|\log r|^C}+ \frac{Cc_1}{|\log r|^C}
\Big(k\int_r^{R}
    |\log s|^{C-1-\mu}\frac{ds}{s}+
\int_r^{R}
    |\log s|^{C-1}\sqrt{s}ds\Big)\leq\\
\frac{c_2}{|\log r|^C}+\frac{kc_2}{|\log r|^C}\int_r^{R}
    |\log s|^{C-1-\mu}\frac{ds}{s}\leq \frac{c_2}{|\log r|^C}+\\
k\frac{c_2}{(C-\mu)|\log r|^C}(|\log r|^{C-\mu}-|\log R|^{C-\mu})
\leq c_2(|\log r|^{-C}+k|\log r|^{-\mu}),
\end{split}\end{equation*}
and the lemma follows.
\end{proof}

\paragraph{Proof of Lemma \ref{frank17}.}\label{17}
 Set
$$
\xi=y-z, \quad u=x-y.
$$
Without loss of generality we can assume that
$$
u=w\mathbf e_1, \quad w=|x-y|, \quad \mathbf e_1=(1,0,0).
$$
We thus get
$$
I=\int_{\mathbb R^3}\frac{|\ln|w\mathbf e_1+\xi||^{\mu/2}
|\ln|\xi||^{\mu/2}}{|w\mathbf e_1+\xi|^2 |\xi|^2} d\xi.
$$
Next set $\zeta=w^{-1}\xi$. It follows that
\begin{equation}\label{frank18}
I=\frac{1}{w}\int_{\mathbb R^3}\frac{|\ln w+\ln|\mathbf
e_1+\zeta||^{\mu/2} |\ln w+\ln|\zeta||^{\mu/2}}{|\mathbf
e_1+\zeta|^2 |\zeta|^2} d\zeta.
\end{equation}
Next notice that
$$
|\ln w+\ln|\mathbf e_1+\zeta||^{\mu/2}\leq c|\ln
w|^{\mu/2}+c|\ln|\mathbf e_1+\zeta||^{\mu/2},
$$
$$
|\ln w+\ln|\zeta||^{\mu/2}\leq c|\ln
w|^{\mu/2}+c|\ln|\zeta||^{\mu/2},
$$
which yields
$$
|\ln w+\ln|\mathbf e_1+\zeta||^{\mu/2} |\ln
w+\ln|\zeta||^{\mu/2}\leq c |\ln w|^{\mu} +c|\ln|\mathbf
e_1+\zeta||^{\mu}+c|\ln|\zeta||^{\mu}
$$
Substituting this inequality into \eqref{frank18} we obtain
\begin{equation}\label{frank19}\begin{split}
 I\leq c \frac{|\ln w|^\mu}{w}\int_{\mathbb R^3}
 \frac{d\zeta}{|\zeta+\mathbf e_1|^2|\zeta|^2}+\\
  c \frac{1}{w}\int_{\mathbb R^3}\frac{|\ln|\mathbf e_1+\zeta||^{\mu}+
  c|\ln|\zeta||^{\mu}}{|\zeta+\mathbf e_1|^2|\zeta|^2}d\zeta\leq c\frac{|\ln w|^\mu+1}{w},
\end{split}\end{equation}
and the lemma follows.

\end{appendix}

%\begin{acknowledgements}
%\end{acknowledgements}




\end{document}